\documentclass[12pt]{article}

\usepackage{amsmath,amssymb,latexsym,pstricks,a4wide,mathrsfs,comment,amsthm,graphicx,tikz}
\usepackage{bibspacing}
\usepackage{geometry}
\begin{document}

\hyphenation{mon-oid mon-oids}
\newcommand{\nc}{\newcommand}
\nc{\rnc}{\renewcommand}

\nc{\partn}[4]{\left( \begin{array}{c|c} 
#1 \ & \ #3 \ \ \\ \cline{2-2}
#2 \ & \ #4 \ \
\end{array} \!\!\! \right)}

\nc{\PXpartn}{\partn{A_i}{B_i}{C_j}{D_k}_{i\in I,\;\! j\in J,\;\! k\in K}}

\nc{\partnL}[2]{(#1|#2)}
\nc{\partnR}[2]{(#1|#2)^*}

\nc{\LXpartn}{\partnL{A_x}{B_x}}
\nc{\LXpartnt}{\partnL{C_x}{D_x}}
\nc{\LXpartnth}{\partnL{E_x}{F_x}}
\nc{\RXpartn}{\partnR{A_x}{B_x}}
\nc{\RXpartnt}{\partnR{C_x}{D_x}}
\nc{\RXpartnth}{\partnR{E_x}{F_x}}

\nc{\sm}{\setminus}
\nc{\De}{\Delta}
\nc{\sups}{\supseteq}
\nc{\G}{\mathscr G}
\rnc{\H}{\mathscr H}
\nc{\an}{\aleph_0}
\nc{\ao}{\aleph_1}
\nc{\N}{\mathbb N}
\nc{\EX}{\E_X}
\nc{\bY}{{\bf Y}}
\nc{\bZ}{{\bf Z}}
\nc{\bW}{{\bf W}}
\nc{\SR}{\operatorname{SR}}
\nc{\alb}{\al^\vee}
\nc{\beb}{\be^\wedge}
\nc{\AND}{\quad\text{and}\quad}
\rnc{\L}{\mathcal L}
\nc{\PXfin}{\P_X^{\fin}}
\nc{\SXfin}{\S_X^{\fin}}

\nc{\J}{\mathcal J}
\nc{\pre}{\preceq}
\nc{\col}{\operatorname{col}}
\nc{\cocol}{\operatorname{cocol}}
\nc{\defect}{\operatorname{def}}
\nc{\codef}{\operatorname{codef}}
\nc{\sing}{\operatorname{sing}}
\nc{\cosing}{\operatorname{cosing}}
\nc{\sh}{\operatorname{sh}}
\nc{\fin}{\operatorname{fin}}

\nc{\bbZ}{\mathbb Z}
\nc{\bbN}{\mathbb N}
\nc{\comma}{,\ }

\rnc{\bar}{\overline}

\nc{\codom}{\operatorname{codom}}
\nc{\coker}{\operatorname{coker}}
\nc{\idrank}{\operatorname{id-rank}}

\nc{\PnSn}{\P_n\setminus\S_n}
\nc{\InSn}{\I_n\setminus\S_n}

\nc{\R}{\mathcal R}

\nc{\suv}[1]{\pscircle*(#1,2){1.5pt}}
\nc{\slv}[1]{\pscircle*(#1,0){1.5pt}}

\nc{\COMMA}{,\quad}
\rnc{\S}{\mathcal S}

\nc{\T}{\mathcal T} 
\nc{\rank}{\operatorname{rank}}

\rnc{\hat}{\widehat}
\nc{\I}{\mathcal I}
\nc{\C}{\mathcal C}
\nc{\B}{\mathcal B}
\nc{\E}{\mathcal E}
\rnc{\P}{\mathcal P}
\nc{\sub}{\subseteq}
\nc{\la}{\langle}
\nc{\ra}{\rangle}
\nc{\mt}{\mapsto}
\nc{\dom}{\mathrm{dom}}
\nc{\im}{\mathrm{im}}
\nc{\id}{\mathrm{id}}
\nc{\ve}{\varepsilon}
\nc{\al}{\alpha}
\nc{\be}{\beta}
\nc{\ga}{\gamma}
\nc{\Ga}{\Gamma}
\nc{\de}{\delta}
\nc{\ka}{\kappa}
\nc{\lam}{\lambda}
\nc{\Lam}{\Lambda}
\nc{\si}{\sigma}
\nc{\vs}{\varsigma}
\nc{\Si}{\Sigma}

\nc{\sC}{\mathscr C}

\nc{\bit}{\vspace{-3 truemm}\begin{itemize}}
\nc{\eit}{\end{itemize}\vspace{-3 truemm}}
\nc{\eitres}{\end{itemize}}

\nc{\set}[2]{\{ {#1} : {#2} \}} 
\nc{\bigset}[2]{\big\{ {#1} : {#2} \big\}} 

\nc{\bE}{{\bf E}}
\nc{\bG}{{\bf G}}

\newcommand{\pf}{\noindent{\bf Proof}\,\,  }
\newcommand{\epf}{\hfill$\Box$\bigskip}
\newcommand{\epfres}{\hfill$\Box$}
\newcommand{\epfreseq}{\tag*{$\Box$}}


\nc{\uv}[1]{\fill (#1,2)circle(.125);}
\nc{\lv}[1]{\fill (#1,0)circle(.125);}
\nc{\stline}[2]{\draw (#1,2)--(#2,0);}
\nc{\uline}[2]{\draw (#1,2)--(#2,2);}
\nc{\lline}[2]{\draw (#1,0)--(#2,0);}
\nc{\udotted}[2]{\draw [dotted] (#1,2)--(#2,2);}
\nc{\ldotted}[2]{\draw [dotted] (#1,0)--(#2,0);}

\nc{\uvc}[2]{\fill [#2] (#1,2)circle(.17);}
\nc{\lvc}[2]{\fill [#2] (#1,0)circle(.15);}

\nc{\uvs}[3]{\fill (#1+#2,2+#3)circle(.125);}
\nc{\lvs}[3]{\fill (#1+#2,0+#3)circle(.125);}
\nc{\uvss}[3]{\fill (#1+#2,2+#3)circle(.1);}
\nc{\lvss}[3]{\fill (#1+#2,0+#3)circle(.1);}
\nc{\stlines}[4]{\draw (#1+#3,2+#4)--(#2+#3,0+#4);}
\nc{\ulines}[4]{\draw (#1+#3,2+#4)--(#2+#3,2+#4);}
\nc{\llines}[4]{\draw (#1+#3,0+#4)--(#2+#3,0+#4);}

\nc{\muv}[1]{\fill (#1,2)circle(.055);}
\nc{\mlv}[1]{\fill (#1,0)circle(.055);}
\nc{\muvs}[3]{\fill (#1+#2,2+#3)circle(.055);}
\nc{\mlvs}[3]{\fill (#1+#2,0+#3)circle(.055);}

\nc{\lblock}[2]{\draw[fill=blue!20] (#1,0)--(#1,.2)--(#2,.2)--(#2,0);}
\nc{\ublock}[2]{\draw[fill=blue!20] (#1,2)--(#1,1.8)--(#2,1.8)--(#2,2);}
\nc{\block}[4]{\draw[fill=blue!20] (#1,2)--(#2,2)--(#4,0)--(#3,0)--(#1,2);}
\nc{\lblocks}[3]{\draw[fill=blue!20] (#1,0+#3)--(#1,.2+#3)--(#2,.2+#3)--(#2,0+#3);}
\nc{\ublocks}[3]{\draw[fill=blue!20] (#1,2+#3)--(#1,1.8+#3)--(#2,1.8+#3)--(#2,2+#3);}
\nc{\blocks}[5]{\draw[fill=blue!20] (#1,2+#5)--(#2,2+#5)--(#4,0+#5)--(#3,0+#5)--(#1,2+#5);}
\nc{\udotteds}[3]{\draw [dotted] (#1,2+#3)--(#2,2+#3);}
\nc{\ldotteds}[3]{\draw [dotted] (#1,0+#3)--(#2,0+#3);}


\nc{\permblock}[4]{\draw[fill=blue!20] (#1,2)--(#2,2) to [out=270,in=90] (#4,0)--(#3,0) to [out=90,in=270] (#1,2);}
\nc{\bpermblock}[4]{\draw[fill=blue!20] (#1,6)--(#2,6) to [out=270,in=90] (#4,0)--(#3,0) to [out=90,in=270] (#1,6);}
\nc{\bpermblockd}[4]{\draw[dotted,line width=0.3pt] (#1,6)--(#2,6) to [out=270,in=90] (#4,0)--(#3,0) to [out=90,in=270] (#1,6);}
\nc{\bpermblocka}[8]{\draw[fill=blue!20] (#1,6)--(#2,6) to [out=270+#5,in=90+#6] (#4,0)--(#3,0) to [out=90+#7,in=270+#8] (#1,6);}
\nc{\permblocks}[5]{\draw[fill=blue!20] (#1,2+#5)--(#2,2+#5) to [out=270,in=90] (#4,0+#5)--(#3,0+#5) to [out=90,in=270] (#1,2+#5);}

\nc{\permline}[2]{\muv{#1}\mlv{#2}\draw (#1,2) to [out=270,in=90] (#2,0);}
\nc{\permlines}[3]{\muvs{#1}{0}{#3}\mlvs{#2}{0}{#3}\draw (#1,2+#3) to [out=270,in=90] (#2,0+#3);}
\nc{\permlinel}[2]{\muv{#1}\mlv{#2}\draw (#1,2) to [out=290,in=110] (#2,0);}
\nc{\lpermlinel}[2]{\muv{#1}\mlv{#2}\draw (#1,2) to [out=250,in=70] (#2,0);}
\nc{\lpermlinels}[3]{\muvs{#1}{0}{#3}\mlvs{#2}{0}{#3}\draw (#1,2+#3) to [out=250,in=70] (#2,0+#3);}
\nc{\permlinels}[3]{\muvs{#1}{0}{#3}\mlvs{#2}{0}{#3}\draw (#1,2+#3) to [out=290,in=110] (#2,0+#3);}
\nc{\permlinells}[3]{\muvs{#1}{0}{#3}\mlvs{#2}{0}{#3}\draw (#1,2+#3) to [out=290,in=120] (#2,0+#3);}
\nc{\permlinellls}[3]{\muvs{#1}{0}{#3}\mlvs{#2}{0}{#3}\draw (#1,2+#3) to [out=250,in=60] (#2,0+#3);}

\newtheorem{thm}[equation]{Theorem}
\newtheorem{lemma}[equation]{Lemma}
\newtheorem{cor}[equation]{Corollary}
\newtheorem{prop}[equation]{Proposition}
\newtheorem{defn}[equation]{Definition}
\newtheorem{eg}[equation]{Example}
\newtheorem{ass}[equation]{Assumption}

\theoremstyle{definition}
\newtheorem{rem}[equation]{Remark}

\title{Infinite partition monoids}
\author{
James East\\
{\footnotesize \emph{School of Computing, Engineering and Mathematics}}\\
{\footnotesize \emph{University of Western Sydney}}\\
{\footnotesize \emph{Locked Bag 1797, Penrith NSW 2751, Australia}}\\
{\footnotesize {\tt J.East\,@\,uws.edu.au}}
}

\date{March 13, 2014}

\maketitle
~\vspace{-1.5 cm}

\begin{abstract}
Let $\P_X$ and $\S_X$ be the partition monoid and symmetric group on an infinite set~$X$.  We show that $\P_X$ may be generated by $\S_X$ together with two (but no fewer) additional partitions, and we classify the pairs $\al,\be\in\P_X$ for which $\P_X$ is generated by $\S_X\cup\{\al,\be\}$.  We also show that $\P_X$ may be generated by the set $\EX$ of all idempotent partitions together with two (but no fewer) additional partitions.  In fact, $\P_X$ is generated by $\EX\cup\{\al,\be\}$ if and only if it is generated by $\EX\cup\S_X\cup\{\al,\be\}$.  We also classify the pairs $\al,\be\in\P_X$ for which $\P_X$ is generated by $\E_X\cup\{\al,\be\}$.  Among other results, we show that any countable subset of $\P_X$ is contained in a $4$-generated subsemigroup of $\P_X$, and that the length function on $\P_X$ is bounded with respect to any generating set.



\end{abstract}



\section{Introduction}

Diagram algebras have been the focus of intense study since the introduction of the Brauer algebras \cite{Bra} in 1937 and, subsequently, the Temperley-Lieb algebras \cite{GW} and Jones algebras~\cite{Jon}.  The \emph{partition algebras}, originally introduced in the context of statistical mechanics~\cite{Mar}, contain all of the above diagram algebras and so provide a unified framework in which to study diagram algebras more generally.  Partition algebras may be thought of as twisted semigroup algebras of \emph{partition monoids}, and many properties of the partition algebras may be deduced from corresponding properties of the associated monoids \cite{JEgrpm,JEpnsn,HR,Wilcox}.  Recent studies have also recognised partition monoids and some of their submonoids as key objects in the pseudovarieties of finite aperiodic monoids and semigroups with involution \cite{Aui2012,Aui2013,ADV}.

Partition monoids were originally defined as finite structures, but the definitions work equally well in the infinite case.  Although most of the study of partition monoids so far has focused on the finite case, there have been a number of recent works on infinite partition monoids; for example, Green's relations were characterized in \cite{FitLau}, and the idempotent generated subsemigroups were described in \cite{EF}.  The purpose of this article is to continue the study of infinite partition monoids, and we investigate a number of problems inspired by analogous considerations in infinite transformation semigroup theory.

As noted in \cite{JEgrpm,EF}, the partition monoids contain a number of important transformation semigroups as submonoids, including the symmetric groups, the full transformation semigroups, and the symmetric and dual symmetric inverse monoids; see \cite{FL,GMbook,Hig,How,Law,Lip,Mal} for background on these subsemigroups.  Many studies of infinite transformation semigroups have concentrated on features concerning generation.  It seems that the earliest result in this direction goes back to 1935, when Sierpi\'nski \cite{Sie} showed that for any infinite set $X$ and for any countable collection $\al_1,\al_2,\ldots$ of functions $X\to X$, it is possible to find functions $\be,\ga:X\to X$ for which each of $\al_1,\al_2,\ldots$ can be obtained by composing $\be$ and $\ga$ in some order a certain number of times.  In modern language, this result says that any countable subset of the full transformation semigroup $\T_X$ is contained in a two-generator subsemigroup, or that the \emph{Sirepi\'nski rank} of infinite $\T_X$ is equal to $2$.  (The Sierpi\'nski rank of a semigroup $S$ is the minimal value of $n$ such that any countable subset of $S$ is contained in an $n$-generator subsemigroup of $S$, if such an $n$ exists, or $\infty$ otherwise.)  Similar results exist for various other transformation semigroups \cite{JEifims,HHMR,HP,MPsurj}; see also \cite{MPsie} for a recent survey.

The notion of Sierpi\'nski rank is intimately connected to the idea of \emph{relative rank}.  
The relative rank of a semigroup $S$ modulo a subset $T\sub S$ is defined to be the least cardinality of a subset~$U$ of $S$ for which $S$ is equal to $\la T\cup U\ra$, the semigroup generated by $T\cup U$.  In the seminal paper on this subject \cite{HHR} (see also \cite{HHR2}), it was shown that an infinite full transformation semigroup $\T_X$ has relative rank $2$ modulo either the symmetric group $\S_X$ or the set $E(\T_X)$ of all idempotents in $\T_X$.  In that paper, the pairs of transformations that, together with $\S_X$ (in the case of $|X|$ being a regular cardinal---see \cite{EMP} for the singular case) or $E(\T_X)$ (for any infinite set $X$), generate all of $\T_X$ were characterized.  Again, these results have led to similar studies of other transformation semigroups \cite{AM,JEifims,HHMR,HMMR,HMR}.

Another closely related concept is the so-called \emph{semigroup Bergman property}; a semigroup has this property if the length function for the semigroup is bounded with respect to any generating set (the bound may be different for different generating sets).  The property is so named because of the seminal paper of Bergman \cite{Berg}, in which it was shown that the infinite symmetric groups have this property;  in fact, Bergman showed that infinite symmetric groups have the corresponding property with respect to \emph{group} generating sets, and the semigroup analogue was proved in \cite{MMR}.  Further studies have investigated the semigroup Bergman property in the context of other transformation semigroups \cite{JEifims,MMR,MPsurj}.  

The goal of the present article is to investigate problems such as those above in the context of infinite partition monoids.
The article is organised as follows.  In Section \ref{sect:preliminaries}, we define the partition monoids $\P_X$ and outline some of their basic properties.  In Section~\ref{sect:rankPXSX}, we show that~$\P_X$ has relative rank $2$ modulo the symmetric group $\S_X$ (Theorem \ref{cor4}) and then, in Section~\ref{sect:genpairsPXSX}, we characterize the pairs $\al,\be\in\P_X$ for which $\P_X$ is generated by $\S_X\cup\{\al,\be\}$.  This characterization depends crucially on the nature of the cardinal $|X|$; we have three separate characterizations, according to whether $|X|$ is countable (Theorem \ref{thm_countable}), or regular but uncountable (Theorem \ref{thm_regular}), or singular (Theorem \ref{thm_singular}).  In Section~\ref{sect:rankPXEX}, we show that the relative rank of $\P_X$ modulo the set $\E_X$ of all idempotent partitions is also equal to $2$ (Theorem~\ref{rankPXEX}); in fact, the relative rank of $\P_X$ modulo $\E_X\cup\S_X$ is equal to $2$ as well.  Then, in Section~\ref{sect:genpairsPXEX}, we show that for any $\al,\be\in\P_X$, $\P_X$ is generated by $\E_X\cup\{\al,\be\}$ if and only if it is generated by $\E_X\cup\S_X\cup\{\al,\be\}$, and we characterize all such pairs $\al,\be$ (Theorem \ref{EXgenthm}).  The characterization in this case does not depend on the cardinality of $X$, but relies crucially on results of~\cite{EF} describing the semigroups $\la \E_X\ra$ and $\la \E_X\cup\S_X\ra$.  Finally, in Section \ref{sect:sierpinski}, we apply the above results to show that $\P_X$ has Sierpi\'nski rank at most~$4$ (Theorem \ref{PXsierpinski}), and also satisfies the semigroup Bergman property (Theorem \ref{PXbergman}).

All functions will be written to the right of their arguments, and functions will be composed from left to right.  We write $A=B\sqcup C$ to indicate that $A$ is the disjoint union of $B$ and~$C$.  We write $\N$ for the set of natural numbers $\{1,2,3,\ldots\}$.  Throughout, a statement such as ``Let~$Y=\set{y_i}{i\in I}$'' should be read as ``Let $Y=\set{y_i}{i\in I}$ and assume the map~$I\to Y:i\mt y_i$ is a bijection''.  We assume the Axiom of Choice throughout.  If $X$ is an infinite set, we will say a family $(X_i)_{i\in I}$ of subsets of $X$ is a \emph{moiety} of $X$ if $X=\bigsqcup_{i\in I}X_i$ and $|X_i|=|X|$ for all $i\in I$.  A cardinal $\mu$ is \emph{singular} if there exists a set $X$ such that $X=\bigcup_{i\in I}X_i$, where $|I|<\mu$ and $|X_i|<\mu$ for each~$i\in I$, but $|X|=\mu$; otherwise, $\mu$ is \emph{regular}.  The only finite regular cardinals are $0$, $1$ and $2$.  The smallest infinite singular cardinal is $\aleph_\omega=\aleph_0+\aleph_1+\aleph_2+\cdots$.  See \cite{Jech} for more details on singular and regular cardinals.

\section{Preliminaries}\label{sect:preliminaries}

In this section, we recall the definition of the partition monoids $\P_X$, and revise some of their basic properties.  We also introduce two submonoids, $\L_X$ and $\R_X$, which will play a crucial role throughout our investigations, and we define a number of parameters associated to a partition that will allow for convenient statements of our results.

Let $X$ be a set, and $X'$ a disjoint set in one-one correspondence with $X$ via a mapping $X\to X':x\mt x'$.  If $A\sub X$ we will write $A'=\set{a'}{a\in A}$. A \emph{partition on $X$} is a collection of pairwise disjoint nonempty subsets of $X\cup X'$ whose union is $X\cup X'$; these subsets are called the \emph{blocks} of the partition.  The \emph{partition monoid} on $X$ is the set $\P_X$ of all partitions on $X$, with a natural associative binary operation defined below. A block $A$ of a partition $\al\in\P_X$ is said to be a \emph{transversal block} if $A\cap X\not=\emptyset\not=A\cap X'$, or otherwise an \emph{upper} (respectively, \emph{lower}) \emph{nontransversal block} if~$A\cap X'=\emptyset$ (respectively, $A\cap X=\emptyset$).  If $\al\in \P_X$, we will write
\[
\al = \partn{A_i}{B_i}{C_j}{D_k}_{i\in I,\;\! j\in J,\;\! k\in K}
\]
to indicate that $\al$ has transversal blocks $A_i\cup B_i'$ ($i\in I$), upper nontransversal blocks $C_j$ (${j\in J}$), and lower nontransversal blocks $D_k'$ ($k\in K$).  The indexing sets $I,J,K$ will sometimes be implied rather than explicit, for brevity; if they are distinct, they will generally be assumed to be disjoint.  Sometimes we will use slight variants of this notation, but it should always be clear what is meant.

\begin{figure}[ht]
   \begin{center}
\begin{tikzpicture}[xscale=.6,yscale=0.6]
  \uv0
  \uv1
  \uv2
  \uv3
  \uv4
  \uv5
  \lv0
  \lv1
  \lv2
  \lv3
  \lv4
  \lv5
  \draw(0,2.1) node [above] {{\tiny $1$}};
  \draw(1,2.1) node [above] {{\tiny $2$}};
  \draw(2,2.1) node [above] {{\tiny $3$}};
  \draw(3,2.1) node [above] {{\tiny $4$}};
  \draw(4,2.1) node [above] {{\tiny $5$}};
  \draw(5,2.1) node [above] {{\tiny $6$}};
  \draw(0.08,-0.1) node [below] {{\tiny $1'$}};
  \draw(1.08,-0.1) node [below] {{\tiny $2'$}};
  \draw(2.08,-0.1) node [below] {{\tiny $3'$}};
  \draw(3.08,-0.1) node [below] {{\tiny $4'$}};
  \draw(4.08,-0.1) node [below] {{\tiny $5'$}};
  \draw(5.08,-0.1) node [below] {{\tiny $6'$}};
  \draw (0,2) .. controls (0,1.2) and (2,1.2) .. (2,2);
  \stline23
  \draw (1,2) .. controls (1,1.2) and (3,1.2) .. (3,2);
  \draw (4,2) .. controls (4,1.5) and (5,1.5) .. (5,2);
  \stline55
  \draw (0,0) .. controls (0,1) and (5,1) .. (5,0);
  \draw (1,0) .. controls (1,0.5) and (2,0.5) .. (2,0);
	\end{tikzpicture}
      \caption{A graphical representation of a partition.}
     \label{fig:partitionfromP6}
   \end{center}
 \end{figure}
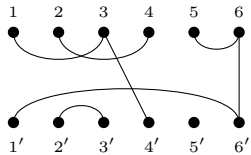

A partition may be represented as a graph on the vertex set $X\cup X'$; edges are included so that the connected components of the graph correspond to the blocks of the partition.  Of course such a graphical representation is not unique, but we regard two such graphs as equivalent if they have the same connected components.  We will also generally identify a partition with any graph representing it.  We think of the vertices from~$X$ (respectively, $X'$) as being the \emph{upper vertices} (respectively, \emph{lower vertices}), explaining our use of these words in relation to the nontransversal blocks.  An example is given in Figure~\ref{fig:partitionfromP6} for the partition $\al=\big\{ \{1,3,4'\},\{2,4\},\{5,6,1',6'\},\{2',3'\},\{5'\}\big\}\in\P_X$, where~$X=\{1,2,3,4,5,6\}$.  Although it is traditional to draw vertex $x$ directly above vertex $x'$, especially in the case of finite $X$, this is not necessary; indeed, we will often be forced to abandon this tradition.  It will also be convenient to sometimes identify a partition $\al\in\P_X$ with its corresponding equivalence relation on $X\cup X'$, and write $(x,y)\in\al$ to indicate that $x,y\in X\cup X'$ belong to the same block of $\al$.

\begin{figure}[ht]
   \begin{center}
\begin{tikzpicture}[xscale=.6,yscale=0.6]
  \uvs004
  \uvs104
  \uvs204
  \uvs304
  \uvs404
  \lvs004
  \lvs104
  \lvs204
  \lvs304
  \lvs404
  \draw (0,2+4) .. controls (0,1.5+4) and (1,1.5+4) .. (1,2+4);
  \draw (1,2+4) .. controls (1,0.7+4) and (3,1.3+4) .. (3,0+4);
  \draw (4,2+4) .. controls (4,0.5+4) and (1,1.5+4) .. (1,0+4);
  \draw (2,0+4) .. controls (2,0.8+4) and (4,0.8+4) .. (4,0+4);
  \draw (-.2,1+4) node [left] {$\al =$};
  \uvs000
  \uvs100
  \uvs200
  \uvs300
  \uvs400
  \lvs000
  \lvs100
  \lvs200
  \lvs300
  \lvs400
  \draw (0,2) .. controls (0,1.5) and (1,1.5) .. (1,2);
  \draw (2,2) .. controls (2,1.2) and (4,1.2) .. (4,2);
  \stlines1100
  \draw (1,0) .. controls (1,0.5) and (2,0.5) .. (2,0);
  \draw (0,0) .. controls (0,1) and (4,1) .. (4,0);
  \draw (-.2,1) node [left] {$\be =$};
  \draw [->] (5,4.75)--(6,4.25);  
  \draw [->] (5,1.25)--(6,1.75);  
  \uvs073
  \uvs173
  \uvs273
  \uvs373
  \uvs473
  \lvs073
  \lvs173
  \lvs273
  \lvs373
  \lvs473
  \draw (0+7,2+3) .. controls (0+7,1.5+3) and (1+7,1.5+3) .. (1+7,2+3);
  \draw (1+7,2+3) .. controls (1+7,0.7+3) and (3+7,1.3+3) .. (3+7,0+3);
  \draw (4+7,2+3) .. controls (4+7,0.5+3) and (1+7,1.5+3) .. (1+7,0+3);
  \draw (2+7,0+3) .. controls (2+7,0.8+3) and (4+7,0.8+3) .. (4+7,0+3);
  \uvs071
  \uvs171
  \uvs271
  \uvs371
  \uvs471
  \lvs071
  \lvs171
  \lvs271
  \lvs371
  \lvs471
  \draw (0+7,2+1) .. controls (0+7,1.5+1) and (1+7,1.5+1) .. (1+7,2+1);
  \draw (2+7,2+1) .. controls (2+7,1.2+1) and (4+7,1.2+1) .. (4+7,2+1);
  \draw (1+7,0+1) .. controls (1+7,0.5+1) and (2+7,0.5+1) .. (2+7,0+1);
  \draw (0+7,0+1) .. controls (0+7,1.0+1) and (4+7,1.0+1) .. (4+7,0+1);
  \stlines1171
  \draw [->] (12,3)--(13,3);  
  \uvs0{14}2
  \uvs1{14}2
  \uvs2{14}2
  \uvs3{14}2
  \uvs4{14}2
  \lvs0{14}2
  \lvs1{14}2
  \lvs2{14}2
  \lvs3{14}2
  \lvs4{14}2
  \draw (0+14,2+2) .. controls (0+14,1.5+2) and (1+14,1.5+2) .. (1+14,2+2);
  \draw (1+14,0+2) .. controls (1+14,0.5+2) and (2+14,0.5+2) .. (2+14,0+2);
  \draw (0+14,0+2) .. controls (0+14,1.0+2) and (4+14,1.0+2) .. (4+14,0+2);
  \draw (4+14,2+2) .. controls (4+14,0.7+2) and (2+14,1.3+2) .. (2+14,0+2);
  \draw (4.2+14,1+2) node [right] {$=\al\be$};
\end{tikzpicture}
    \caption{Two partitions $\al,\be$ (left), their product $\al\be$ (right), and the product graph $\Ga(\al,\be)$ (centre).}
    \label{fig:multinP5}
   \end{center}
 \end{figure}
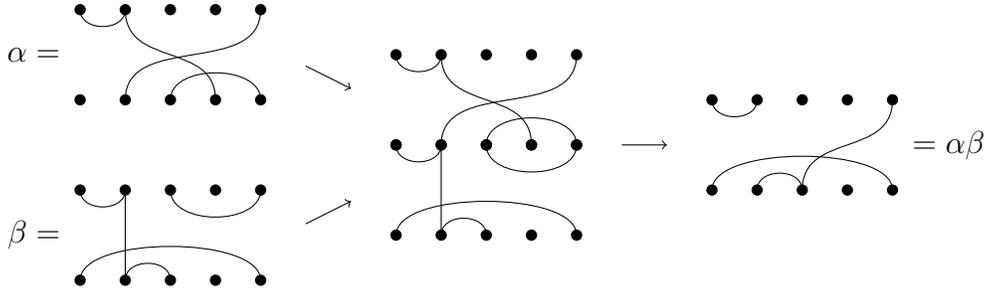

The rule for multiplication of partitions is best described in terms of the graphical representations.  Let $\al,\be\in\P_X$.  Consider now a third set $X''$, disjoint from both $X$ and $X'$, and in bijection with both sets via the maps $X\to X'':x\mt x''$ and $X'\to X'':x'\mt x''$.  Let $\alb$ be the graph obtained from (a graph representing) $\al$ simply by changing the label of each lower vertex $x'$ to $x''$.  Similarly, let $\beb$ be the graph obtained from $\be$ by changing the label of each upper vertex $x$ to $x''$.  Consider now the graph $\Ga(\al,\be)$ on the vertex set~$X\cup X'\cup X''$ obtained by joining $\alb$ and $\beb$ together so that each lower vertex $x''$ of $\alb$ is identified with the corresponding upper vertex $x''$ of $\beb$.  Note that $\Ga(\al,\be)$, which we call the \emph{product graph} of $\al$ and $\be$, may contain multiple edges.  We define $\al\be\in\P_X$ to be the partition that satisfies the property that $x,y\in X\cup X'$ belong to the same block of $\al\be$ if and only if there is a path from $x$ to $y$ in $\Ga(\al,\be)$.  An example calculation (with $X$ finite) is given in Figure~\ref{fig:multinP5}.  (See also~\cite{Mal} for an equivalent formulation of the product; there $\P_X$ was denoted $\mathcal{CS}_X$, and called the \emph{composition semigroup} on $X$.)

This product is easily checked to be associative, and so gives $\P_X$ the structure of a monoid; the identity element is the partition $\bigset{\{x,x'\}}{x\in X}$, which we denote by $1$.  A partition $\al\in\P_X$ is a unit if and only if each block of $\al$ is of the form $\{x,y'\}$ for some $x,y\in X$.  So it is clear that the group of units, which we denote by $\S_X$, is (isomorphic to) the symmetric group on $X$.  So, if $\pi\in\S_X$ and $x\in X$, we will write $x\pi$ for ``the image of $x$ under $\pi$'', by which we mean the unique element of $X$ such that $\{x,(x\pi)'\}$ is a block of~$\pi$.


A crucial aspect of the structure of $\P_X$ is given by the map ${}^*:\P_X\to\P_X:\al\mt\al^*$ where $\al^*$ is the result of ``turning $\al$ upside-down''.  More precisely:
\[
\al=\partn{A_i}{B_i}{C_j}{D_k} \quad\Rightarrow\quad \al^*=\partn{B_i}{A_i}{D_k}{C_j}.
\]
Note that $\pi^*=\pi^{-1}$ if $\pi\in\S_X$.
The next lemma is proved easily, and collects the basic properties of the ${}^*$ map that we will need.  Essentially it states that $\P_X$ is a \emph{regular $*$-semigroup}.

\bigskip
\begin{lemma}\label{*-regular}
Let $\al,\be\in\P_X$.  Then
\begin{equation}\epfreseq
(\al^*)^*=\al\COMMA \al\al^*\al = \al\COMMA \al^*\al\al^*=\al^*\COMMA (\al\be)^*=\be^*\al^*.
\end{equation}
\end{lemma}

Among other things, these properties mean that the map $\al\mt\al^*$ is an anti-isomorphism of $\P_X$.  This duality will allow us to shorten many proofs.

Next we record some notation and terminology.  With this in mind, let $\al\in\P_X$.  For ${x\in X\cup X'}$, we denote the block of $\al$ containing $x$ by $[x]_\al$.  The \emph{domain} and \emph{codomain} of~$\al$ are defined to be the following subsets of $X$:
\begin{align*}
\dom(\al) &= \bigset{ x\in X } { [x]_\al\cap X'\not=\emptyset}, \\
\codom(\al) &= \bigset{ x\in X } { [x']_\al\cap X\not=\emptyset}.\\
\intertext{We also define the \emph{kernel} and \emph{cokernel} of $\al$ to be the following equivalences on $X$:}
\ker(\al) &= \bigset{(x,y)\in X\times X}{[x]_\al=[y]_\al}, \\
\coker(\al) &= \bigset{(x,y)\in X\times X}{[x']_\al=[y']_\al}.
\end{align*}
Note that $\dom(\al^*)=\codom(\al)$ and $\ker(\al^*)=\coker(\al)$.

\newpage

\begin{lemma}\label{lem1}
Let $\al,\be\in\P_X$.  Then
\bit
	\item[\emph{(\ref{lem1}.1)}] $\dom(\al\be)\sub\dom(\al)$, with equality if $\codom(\al)\sub\dom(\be)$,
	\item[\emph{(\ref{lem1}.2)}] $\codom(\al\be)\sub\codom(\be)$, with equality if $\dom(\be)\sub\codom(\al)$,
	\item[\emph{(\ref{lem1}.3)}] $\ker(\al\be)\sups\ker(\al)$, with equality if $\ker(\be)\sub\coker(\al)$, and
	\item[\emph{(\ref{lem1}.4)}] $\coker(\al\be)\sups\coker(\be)$, with equality if $\coker(\al)\sub\ker(\be)$.
\eit
\end{lemma}

\pf We will only prove (\ref{lem1}.1) and (\ref{lem1}.3), since the others follow by duality.
Clearly $\dom(\al\be)\sub\dom(\al)$.  Suppose $\codom(\al)\sub\dom(\be)$.  Let $x\in\dom(\al)$.  Then $(x,y')\in\al$ for some $y\in\codom(\al)$.  Since $\codom(\al)\sub\dom(\be)$, it follows that $(y,z')\in\be$ for some $z\in\codom(\be)$.  Then $(x,z')\in\al\be$, whence $x\in\dom(\al\be)$, establishing (\ref{lem1}.1).

Clearly $\ker(\al\be)\sups\ker(\al)$.  Suppose $\ker(\be)\sub\coker(\al)$.  Let $(x,y)\in\ker(\al\be)$.  If one of $x$ or $y$ belongs to $X\sm\dom(\al)$, then so too does the other, and $(x,y)\in\ker(\al)$.  So suppose $x,y\in\dom(\al)$.  Then $(x,a'),(y,b')\in\al$ for some $a,b\in\codom(\al)$.  Since $(x,y)\in\ker(\al\be)$, there exist $x_0,x_1,\ldots,x_r\in X$ such that $x_0=a$, $x_r=b$ and $(x_0,x_1)\in\coker(\al)$, $(x_1,x_2)\in\ker(\be),(x_2,x_3)\in\coker(\al)$, and so on.  But, since $\ker(\be)\sub\coker(\al)$, it follows that $(x_0,x_1),(x_1,x_2),\ldots,(x_{r-1},x_r)\in\coker(\al)$.  This then implies that $(a,b)\in\coker(\al)$, and $(x,y)\in\ker(\al)$.  This completes the proof of (\ref{lem1}.3).  \epf

We now define two submonoids of $\P_X$ that will play a crucial role in what follows.  Denote by $\De=\bigset{(x,x)}{x\in X}$ the trivial equivalence (that is, the equality relation).  Let
\begin{align*}
\L_X &= \set{\al\in\P_X}{\dom(\al)=X,\ \ker(\al)=\De},\ \text{and} \\
\R_X &= \set{\al\in\P_X}{\codom(\al)=X,\ \coker(\al)=\De}.
\end{align*}
Note that $\L_X^*=\R_X$ and $\R_X^*=\L_X$, and that $\L_X\cap\R_X=\S_X$.

\bigskip
\begin{lemma}\label{lem2}
The sets $\L_X$ and $\R_X$ are submonoids of $\P_X$.  Further, $\P_X\sm\L_X$ is a right ideal of $\P_X$, and $\P_X\sm\R_X$ is a left ideal.
\end{lemma}

\pf We will prove the statements concerning $\L_X$, and those concerning $\R_X$ will follow by duality.  Let $\al,\be\in\L_X$.  Then $\dom(\al\be)=\dom(\al)=X$ and $\ker(\al\be)=\ker(\al)=\De$ by (\ref{lem1}.1) and (\ref{lem1}.3), respectively, so that $\al\be\in\L_X$.  

Next, let $\al\in\P_X\sm\L_X$ and $\be\in\P_X$.  
%
%
If $\dom(\al)\not=X$, then $\dom(\al\be)\sub\dom(\al)\not=X$, so that $\dom(\al\be)\not=X$, and $\al\be\in\P_X\sm\L_X$.  If $\ker(\al)\not=\De$, then we similarly obtain $\al\be\in\P_X\sm\L_X$.~\epf

\begin{rem}
As noted in \cite{JEgrpm,EF}, the submonoids $\set{\al\in\P_X}{\dom(\al)=\codom(\al)=X}$ and $\set{\al\in\P_X}{\ker(\al)=\coker(\al)=\De}$ are isomorphic to the symmetric inverse semigroup and dual symmetric inverse semigroup on $X$, respectively.
\end{rem}

\begin{figure}[ht]
   \begin{center}
\begin{tikzpicture}[xscale=.8,yscale=0.8]
  \muvs002
  \blocks000{1.5}2
  \muv0
  \muv{0.25}
  \muv{0.5}
  \block000{.5}
  \block{.25}{.25}1{1.5}
  \block{.5}{.5}{2}{2.5}
  \udotteds{.95}{1.5}{-.1}
  \ldotted{2.75}4
  \draw(0,4)node[above]{\small $x$};
	 \draw(2,0.1)node[below]{{\small $\underbrace{\phantom{...............................}}_{E_x=\bigcup_{y\in A_x}C_y}$}};   %
  \lblocks{6}{7.5}2
  \muv6
  \muv{6.25}
  \muv{6.5}
  \block666{6.5}
  \block{6.25}{6.25}7{7.5}
  \block{6.5}{6.5}{8}{8.5}
  \udotteds{6.95}{7.5}{-.1}
  \ldotted{8.75}{10}
	 \draw(8,0.1)node[below]{{\small $\underbrace{\phantom{...............................}}_{F_i=\bigcup_{y\in B_i}C_y}$}};   	 \draw(6.75,2.1)node[above]{{\small $\overbrace{\phantom{............}}^{B_i}$}};   %
  \draw[|-]  (0-3,4)--(0-3,0);
  \draw[|-|] (0-3,2)--(0-3,0);
  \draw(0-3,3)node[left]{{\small $\al$}};
  \draw(0-3,1)node[left]{{\small $\be$}};
	\end{tikzpicture}
    \caption{The product $\al\be=\partnL{E_x}{F_i,D_j}$ of two elements $\al=\partnL{A_x}{B_i}$ and $\be=\partnL{C_x}{D_j}$ from~$\L_X$, focusing on the blocks $\{x\}\cup E_x'$ (left) and $F_i'$ (right).  See text for further explanation.}
    \label{fig:multinLX}
   \end{center}
 \end{figure}
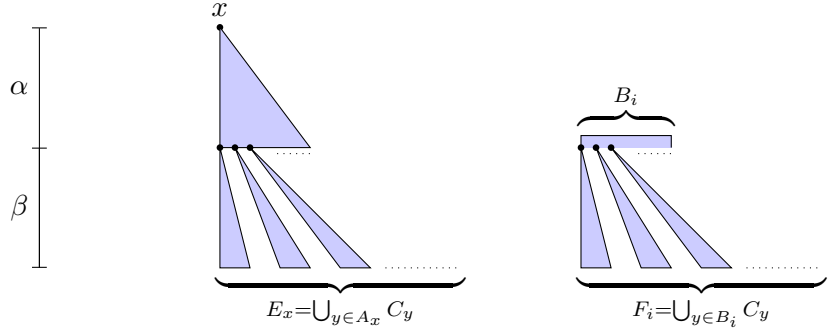

A typical element of $\L_X$ has the form
\[
\partn{x}{A_x}{\emptyset}{B_i}_{x\in X,\;\! i\in I}.
\]
In what follows, we will shorten this to $\partnL{A_x}{B_i}_{x\in X,i\in I}$, or just $\partnL{A_x}{B_i}$.  Accordingly, we will write $\partnR{A_x}{B_i}$ for the partition
\[
\partn{A_x}{x}{B_i}{\emptyset}_{x\in X,\;\! i\in I}
\]
from $\R_X$.  Note that if $\al=\partnL{A_x}{B_i}$ and $\be=\partnL{C_x}{D_j}$, then $\al\be=\partnL{E_x}{F_i,D_j}$, where
$
E_x = \bigcup_{y\in A_x}C_y $ and $ F_i = \bigcup_{y\in B_i}C_y
$
for each $x\in X$ and $i\in I$; see Figure \ref{fig:multinLX}.  A similar rule holds for multiplication in $\R_X$.

We now define a number of parameters associated with a partition.  With this in mind, let $\al\in\P_X$ and write
\[
\al = \partn{A_i}{B_i}{C_j}{D_k}_{i\in I,\;\! j\in J,\;\! k\in K}.
\]
For any cardinal $\mu\leq|X|$, we define
\[
\begin{array}{rclcrcl}
k(\al,\mu) \hspace{-.25cm}&=&\hspace{-.25cm} \# \bigset{i\in I}{|A_i|\geq\mu}, & &
d(\al,\mu) \hspace{-.25cm}&=&\hspace{-.25cm} \# \bigset{j\in J}{|C_j|\geq\mu},\\
k^*(\al,\mu) \hspace{-.25cm}&=&\hspace{-.25cm} \# \bigset{i\in I}{|B_i|\geq\mu}, & &
d^*(\al,\mu) \hspace{-.25cm} &=&\hspace{-.25cm} \# \bigset{k\in K}{|D_k|\geq\mu}.
\end{array}
\]
Note that $k^*(\al,\mu)=k(\al^*,\mu)$ and $d^*(\al,\mu)=d(\al^*,\mu)$.  We also have identities such as $d(\al,\mu)\geq d(\al,\nu)$ if $\mu\leq\nu\leq|X|$.
It will also be convenient to write
\[
d(\al)=d(\al,1)=|J| \AND \ d^*(\al)=d^*(\al,1)=|K|.
\]
The above parameters are natural extensions of those introduced in the context of transformation semigroups in \cite{How66} (see also \cite{EMP,HHR}).  These parameters should not be confused with those introduced in \cite{EF}, such as $\defect(\al)$, $\col(\al)$, etc.

\bigskip
\begin{lemma}\label{lem5.1}
Let $\al,\be\in\P_X$.  Then
\bit
	\item[\emph{(\ref{lem5.1}.1)}] $d(\al)\leq d(\al\be) \leq d(\al)+ d(\be)$, and
  \item[\emph{(\ref{lem5.1}.2)}] $d^*(\be)\leq d^*(\al\be) \leq d^*(\al)+ d^*(\be)$.
\eit
\end{lemma}

\pf We just prove (\ref{lem5.1}.1), since (\ref{lem5.1}.2) will follow by duality.  Let 
\[
\al = \partn{A_i}{B_i}{C_j}{D_k}_{i\in I,\;\! j\in J,\;\! k\in K}
\AND
\be = \partn{E_l}{F_l}{G_m}{H_n}_{l\in L,\;\! m\in M,\;\! n\in N}.
\]
Note that each $C_j$ is an upper nontransversal block of $\al\be$, so $d(\al\be)\geq d(\al)$.  Suppose now that $P$ is an upper nontransversal block of $\al\be$ but that $P\not=C_j$ for any $j\in J$.  Then $P=\bigcup_{i\in I_P}A_i$ for some subset $I_P\sub I$.  Now, $\bigcup_{i\in I_P}B_i$ must have trivial intersection with each of the $E_l$, or else $P$ would be contained in a transversal block of $\al\be$.  But this implies that $\bigcup_{i\in I_P}B_i$ intersects at least one of the $G_m$.  In particular, there are at most $|M|$ such upper nontransversal blocks $P$.  Thus, $d(\al\be)\leq|J|+|M|=d(\al)+d(\be)$. \epf

\begin{rem}
The above-mentioned rule for multiplication in $\L_X$ shows that $d^*(\al\be)=d^*(\al)+d^*(\be)$ if $\al,\be\in\L_X$.  A dual identity holds in $\R_X$.
\end{rem}

For the following lemmas, recall that we count $1$ and $2$ as regular cardinals.

\bigskip
\begin{lemma}\label{k-lemma}
Let $\al,\be\in\L_X$ and $\ga,\de\in\R_X$ and let $\mu\leq|X|$ be any cardinal.  Then
\bit
	\item[\emph{(\ref{k-lemma}.1)}] $k^*(\al,\mu)\leq k^*(\al\be,\mu)$,
	\item[\emph{(\ref{k-lemma}.2)}] $k^*(\al\be,\mu) \leq k^*(\al,\mu)+k^*(\be,\mu)$ if $\mu$ is regular,
	\item[\emph{(\ref{k-lemma}.3)}] $k(\de,\mu)\leq k(\ga\de,\mu)$, and
	\item[\emph{(\ref{k-lemma}.4)}] $k(\ga\de,\mu) \leq k(\ga,\mu)+k(\de,\mu)$ if $\mu$ is regular.
\eit
\end{lemma}

\pf We just prove (\ref{k-lemma}.1) and (\ref{k-lemma}.2), since the others will follow by duality.  Let $\al=\partnL{A_x}{B_i}$ and $\be=\partnL{C_x}{D_j}$.
Then $\al\be=\partnL{E_x}{F_i,D_j}$,
where $E_x = \bigcup_{y\in A_x}C_y$ and $F_i = \bigcup_{y\in B_i}C_y$ for each $x\in X$ and $i\in I$.  Clearly, $|E_x|\geq|A_x|$ for all $x\in X$, so $k^*(\al\be,\mu)\geq k^*(\al,\mu)$, establishing~(\ref{k-lemma}.1).  Next, suppose $\mu$ is regular.  If $|E_x|\geq\mu$ for some $x\in X$, then either (i)~$|A_x|\geq\mu$, or (ii)~$|C_y|\geq\mu$ for some $y\in A_x$.  There are $k^*(\al,\mu)$ values of $x$ that satisfy~(i), and at most $k^*(\be,\mu)$ values of $x$ that satisfy (ii).  Thus, $k^*(\al\be,\mu)\leq k^*(\al,\mu)+k^*(\be,\mu)$, establishing (\ref{k-lemma}.2).~\epf

\newpage

\begin{lemma}\label{d-lemma}
Let $\al,\be\in\L_X$ and $\ga,\de\in\R_X$ and let $\mu\leq|X|$ be any cardinal.  Then
\bit
	\item[\emph{(\ref{d-lemma}.1)}] $d^*(\be,\mu)\leq d^*(\al\be,\mu)$,
	\item[\emph{(\ref{d-lemma}.2)}] $d^*(\al\be,\mu) \leq d^*(\al,\mu)+d^*(\be,\mu)+k^*(\be,\mu)$ if $\mu$ is regular,
	\item[\emph{(\ref{d-lemma}.3)}] $d(\ga,\mu)\leq d(\ga\de,\mu)$, and
	\item[\emph{(\ref{d-lemma}.4)}] $d(\ga\de,\mu) \leq d(\ga,\mu)+d(\de,\mu)+k(\ga,\mu)$ if $\mu$ is regular.
\eit
\end{lemma}

\pf Again, it suffices to prove (\ref{d-lemma}.1) and (\ref{d-lemma}.2).  Write $\al=\partnL{A_x}{B_i}$ and $\be=\partnL{C_x}{D_j}$.  Then $\al\be=\partnL{E_x}{F_i,D_j}$, where $E_x = \bigcup_{y\in A_x}C_y$ and $F_i = \bigcup_{y\in B_i}C_y$ for each $x\in X$ and $i\in I$.  There are $d^*(\be,\mu)$ values of $j\in J$ for which $|D_j|\geq\mu$.  It follows that $d^*(\al\be,\mu)\geq d^*(\be,\mu)$.  Next, suppose $\mu$ is regular, and that $i\in I$ is such that $|F_i|\geq\mu$.  Then either (i)~$|B_i|\geq\mu$, or (ii)~$|C_y|\geq\mu$ for some $y\in B_i$.  There are $d^*(\al,\mu)$ values of $i$ for which (i) holds, and at most $k^*(\be,\mu)$ values of $i$ for which (ii) holds.  Thus, $d^*(\al\be,\mu)\leq d^*(\al,\mu)+d^*(\be,\mu)+k^*(\be,\mu)$.~\epf

The next lemma will be used on a number of occasions.  There is a dual result, but we will not need to state it.

\bigskip
\begin{lemma}\label{lem5file}
Let $\al\in\L_X$ with $d^*(\al)=|X|$, and let $\mu\leq|X|$ be any cardinal.  Then there exists $\be\in\la\S_X,\al\ra\sub\L_X$ such that $d^*(\be,\mu)\geq k^*(\al,\mu)$, $d^*(\be)=|X|$, and $|x|_\be\geq|x|_\al$ for all $x\in X$.
\end{lemma}

\pf Let
$
\al=\LXpartn
$
and put $Y=\set{x\in X}{|A_x|\geq\mu}$, noting that $|Y|=k^*(\al,\mu)$.  We will consider two separate cases.

{\bf Case 1.}  First suppose $|Y|<|X|$.  For each $x\in Y$, choose some $b_x\in B_x$.  Let $\pi\in\S_X$ be any permutation that extends the map $\set{b_x}{x\in Y}\to Y:b_x\mt x$, and put $\be=\al\pi\al$.  Then, for each $x\in Y$, $\left(\bigcup_{y\in B_x}A_{y\pi}\right)'$ is a lower nontransversal block of $\be$, and $\left|\bigcup_{y\in B_x}A_{y\pi}\right|\geq|A_{b_x\pi}|=|A_x|\geq\mu$.  Thus, $d^*(\be,\mu)\geq|Y|$.

{\bf Case 2.}  Now suppose $|Y|=|X|$.  Let $(Y_1,Y_2)$ be a moiety of $Y$, and let $\pi\in\S_X$ be any permutation that extends any bijection $\bigcup_{x\in X}B_x\to Y_1$.  Then for any $x\in X$, $\left(\bigcup_{y\in B_x}A_{y\pi}\right)'$ is a lower nontransversal block of $\be=\al\pi\al$ of size at least $\mu$.  It follows that ${d^*(\be,\mu)=|X|=|Y|}$.

In either case, $[x]_\be=\{x\}\cup\left(\bigcup_{y\in A_x}A_{y\pi}\right)'$, so that $|x|_\be\geq1+|A_x|=|x|_\al$ for all $x\in X$.  And, in either case, $d^*(\be)=|X|$ is a consequence of (\ref{lem5.1}.2). \epf



\section{Relative rank of $\P_X$ modulo $\S_X$}\label{sect:rankPXSX}

Recall that the \emph{relative rank} of a semigroup $S$ with respect to a subset $T$, denoted $\rank(S:T)$, is the minimum cardinality of a subset $U\sub S$ such that $S=\la T\cup U\ra$.  Our goal in this section is to show that $\rank(\P_X:\S_X)=2$; see Theorem \ref{cor4}.

Recall that for $\al\in\P_X$ and $x\in X\cup X'$, we write $[x]_\al$ for the block of $\al$ containing $x$.  We will also write $|x|_\al$ for the cardinality of $[x]_\al$.  The next result shows that $\P_X$ may be generated by $\S_X$ along with just two additional partitions.  See \cite[Theorem 3.3]{HHR} for the corresponding result for infinite transformation semigroups.

\bigskip
\begin{prop}\label{thm3}
Let $\al\in\L_X$ and $\be\in\R_X$ be such that $d^*(\al,|X|)=d(\be,|X|)=|X|$, and $|x|_\al=|x'|_\be=|X|$ for all $x\in X$.  Then $\P_X
=\la\S_X,\al,\be\ra$.
\end{prop}

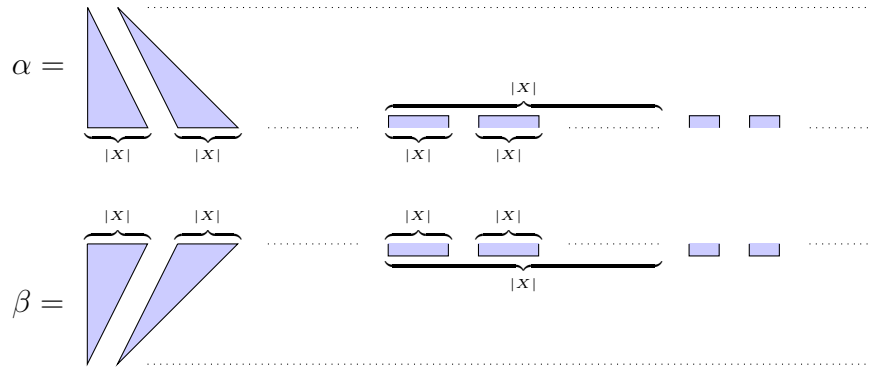
\begin{figure}[ht]
   \begin{center}
\begin{tikzpicture}[xscale=.8,yscale=0.8]
	 \draw(-0.2,1)node[left]{$\al=$}; 
	 \draw(.5,0.2)node[below]{{\tiny $\underbrace{\phantom{...........}}_{|X|}$}}; 
	 \draw(2,0.2)node[below]{{\tiny $\underbrace{\phantom{...........}}_{|X|}$}}; 
	 \draw(5.5,0.2)node[below]{{\tiny $\underbrace{\phantom{...........}}_{|X|}$}}; 
	 \draw(7,0.2)node[below]{{\tiny $\underbrace{\phantom{...........}}_{|X|}$}}; 
	 \draw(7.25,0)node[above]{{\tiny $\overbrace{\phantom{..............................................}}^{|X|}$}}; 
  \block0001
  \block{.5}{.5}{1.5}{2.5}
  \lblock56
  \lblock{6.5}{7.5}
  \lblock{10}{10.5}
  \lblock{11}{11.5}
  \ldotted3{4.5}
  \ldotted{8}{9.5}
  \udotted1{13}
  \ldotted{12}{13}
	 \draw(7,-.5)node[below]{{\tiny $\phantom{|X|}$}}; 
	\end{tikzpicture}
	~ ~ ~
	\begin{tikzpicture}[xscale=.8,yscale=0.8]
		 \draw(-0.2,1)node[left]{$\be=$}; 
	 \draw(.5,1.8)node[above]{{\tiny $\overbrace{\phantom{...........}}^{|X|}$}}; 
	 \draw(2,1.8)node[above]{{\tiny $\overbrace{\phantom{...........}}^{|X|}$}}; 
	 \draw(5.5,1.8)node[above]{{\tiny $\overbrace{\phantom{...........}}^{|X|}$}}; 
	 \draw(7,1.8)node[above]{{\tiny $\overbrace{\phantom{...........}}^{|X|}$}}; 
	 \draw(7.25,2)node[below]{{\tiny $\underbrace{\phantom{..............................................}}_{|X|}$}}; 
  \block0100
  \block{1.5}{2.5}{.5}{.5}
  \ublock56
  \ublock{6.5}{7.5}
  \ublock{10}{10.5}
  \ublock{11}{11.5}
  \udotted3{4.5}
  \udotted{8}{9.5}
  \ldotted1{13}
  \udotted{12}{13}
	\end{tikzpicture}
    \caption{The partitions $\al$ (top) and $\be$ (bottom) from the proof of Proposition~\ref{thm3}.}
    \label{fig:genpair}
   \end{center}
 \end{figure}

\pf Consider an arbitrary partition
\[
\ga = \partn{A_i}{B_i}{C_j}{D_k}_{i\in I,\;\! j\in J,\;\! k\in K}.
\]
We will construct a permutation $\pi\in\S_X$ such that $\ga=\al\pi\be$.  The assumptions on $\al,\be$ allow us to write $\al= \partnL{E_x}{F_x}$ and $\be= \partnR{G_x}{H_x}$, where $|E_x|=|G_x|=|X|$ for all $x\in X$.  See Figure \ref{fig:genpair} for an illustration (the picture shows the basic ``shape'' of $\al$ and $\be$, and is not meant to indicate that $\be=\al^*$).
Let
\[
\begin{array}{rclcrcl}
X_1 \hspace{-.25cm}&=&\hspace{-.25cm} \# \bigset{x\in X}{|F_x|=|X|}, & &
X_3 \hspace{-.25cm}&=&\hspace{-.25cm} \# \bigset{x\in X}{|H_x|=|X|},\\
X_2 \hspace{-.25cm}&=&\hspace{-.25cm} \# \bigset{x\in X}{|F_x|<|X|}, & &
X_4 \hspace{-.25cm}&=&\hspace{-.25cm} \# \bigset{x\in X}{|H_x|<|X|}.
\end{array}
\]
So $X=X_1\sqcup X_2 = X_3\sqcup X_4$.  Note that $|X_1|=d^*(\al,|X|)=|X|$ and, similarly, $|X_3|=|X|$.  We now proceed to construct $\pi\in\S_X$ in stages.

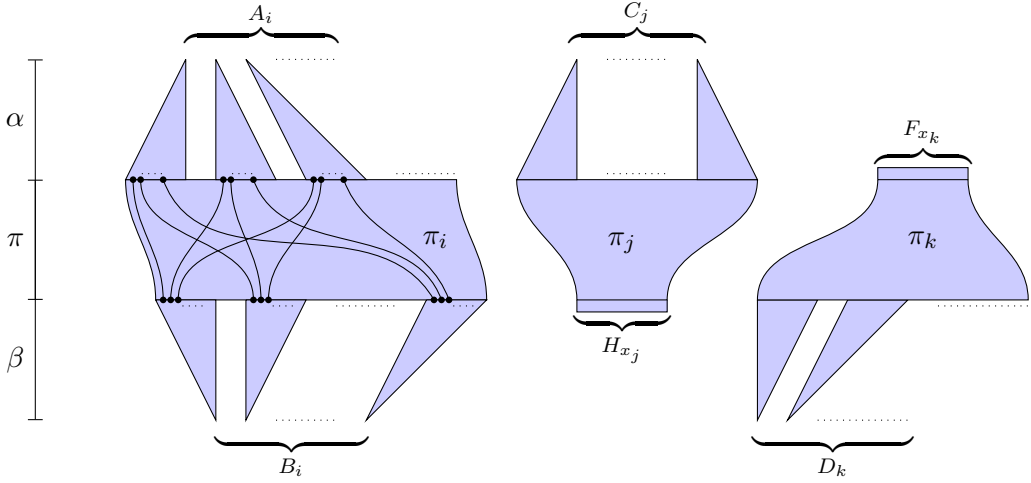
\begin{figure}[ht]
   \begin{center}
\begin{tikzpicture}[xscale=.8,yscale=0.8]
  \permblocks{1.5}{7}{2}{7.5}2
  \blocks{2.5}{2.5}{1.5}{2.5}4
  \blocks33344
  \blocks{3.5}{3.5}{4.5}{5.5}4
  \udotteds454
  \ldotteds67{4.1}
  \draw(3.75,6)node[above]{{\small $\overbrace{\phantom{...................}}^{A_i}$}};
  \block2333
  \block{3.5}{4.5}{3.5}{3.5}
  \block{6.5}{7.5}{5.5}{5.5}
  \draw(4.25,0)node[below]{{\small $\underbrace{\phantom{...................}}_{B_i}$}};
  \permlines{1.75-.125}{2.25-.125}2
  \permlines{2-.125-.125}{3.75-.125}2
  \permlines{1.75+1.5-.125}{2.25+.25-.125-.125}2
  \permlines{2+1.5-.125-.125}{3.75+.25-.125-.125}2
  \permlines{1.75+1.5+1.5-.125}{2.25+.25+.25-.125-.125-.125}2
  \permlines{2+1.5+1.5-.125-.125}{3.75+.25+.25-.125-.125-.125}2
  \permlinells{2.125}{6.75-.125}2
  \permlinels{2+1.5-.125+.25}{6.75+.25-.125-.125}2
  \permlinels{2+1.5+1.5-.125+.25}{6.75+.25+.25-.125-.125-.125}2
  \udotteds56{-.1}
  \ldotted45
  \ldotteds{1.77}{2.15}{4.1}
  \ldotteds{1.77+1.5}{2.15+1.5}{4.1}
  \ldotteds{1.77+3}{2.15+3}{4.1}
  \udotteds{1.77+.675}{2.15+.675}{-.1}
  \udotteds{1.77+.675+1.5}{2.15+.675+1.5}{-.1}
  \udotteds{1.77+.675+4.5}{2.15+.675+4.5}{-.1}
  \ublock9{10.5}
  \permblocks8{12}9{10.5}2
  \blocks99894
  \blocks{11}{11}{11}{12}4
  \udotteds{9.5}{10.5}4
  \ldotteds{9.5}{10.5}{4.1}
  \draw(10,6)node[above]{{\small $\overbrace{\phantom{.................}}^{C_j}$}};
  \draw(9.75,2)node[below]{{\small $\underbrace{\phantom{............}}_{H_{x_j}}$}};
  \draw(9.75,3-.05)node{{\small $\pi_j$}};
  \block{12}{13}{12}{12}
  \block{13.5}{14.5}{12.5}{12.5}
  \ldotted{13}{14.5}
  \udotteds{15}{16.5}{-.1}
  \lblocks{14}{15.5}4
  \permblocks{14}{15.5}{12}{16.5}2
  \draw(14.75,4)node[above]{{\small $\overbrace{\phantom{............}}^{F_{x_k}}$}};
  \draw(13.25,0)node[below]{{\small $\underbrace{\phantom{....................}}_{D_k}$}};
  \draw(14.75,3)node{{\small $\pi_k$}};
  \draw[|-|] (0,6)--(0,0);
  \draw[|-|] (0,4)--(0,2);
  \draw(0,5)node[left]{{\small $\al$}};
  \draw(0,3)node[left]{{\small $\pi$}};
  \draw(0,1)node[left]{{\small $\be$}};
  \draw(6.65,3)node{{\small $\pi_i$}};
	\end{tikzpicture}
    \caption{A schematic diagram of the product $\al\pi\be$, focusing on a transversal block $A_i\cup B_i'$ (left), an upper nontransversal block $C_j$ (middle), and a lower nontransversal block $D_k'$ (right).  See text for further explanation.}
    \label{fig:stage1}
   \end{center}
 \end{figure}

{\bf Stage 1.}  Fix $i\in I$.  For each $x\in A_i$, let $E_x=E_x^1\sqcup E_x^2$ where $|E_x^1|=|B_i|$ and $|E_x^2|=|X|$, and write $E_x^1=\set{e_{xy}}{y\in B_i}$.  For each $y\in B_i$, let $G_y=G_y^1\sqcup G_y^2$ where $|G_y^1|=|A_i|$ and $|G_y^2|=|X|$, and write $G_y^1=\set{g_{xy}}{x\in A_i}$.  Now let
$
\pi_i:\bigcup_{x\in A_i} E_x \to \bigcup_{y\in B_i} G_y
$
be any bijection that extends the map
$
\bigcup_{x\in A_i} E_x^1 \to \bigcup_{y\in B_i} G_y^1: e_{xy}\mt g_{xy}.
$
It is easy to check that if $\pi\in\S_X$ is any permutation that extends $\pi_i$, then $A_i\cup B_i'$ is a block of $\al\pi\be$.  See Figure \ref{fig:stage1} for an illustration.

{\bf Stage 2.}  Let $X_3=X_3^1\sqcup X_3^2$ where $|X_3^1|=|J|$ and $X_3^2\not=\emptyset$, and write $X_3^1=\set{x_j}{j\in J}$.  Now fix $j\in J$.  Choose any bijection $\pi_j:\bigcup_{x\in C_j}E_x\to H_{x_j}$.  Again, it is easy to check that if $\pi\in\S_X$ is any permutation that extends $\pi_j$, then $C_j$ is a block of $\al\pi\be$.  See Figure \ref{fig:stage1}.

{\bf Stage 3.}  Let $X_1=X_1^1\sqcup X_1^2$ where $|X_1^1|=|K|$ and $X_1^2\not=\emptyset$, and write $X_1^1=\set{x_k}{k\in K}$.  Now fix $k\in K$.  Choose any bijection $\pi_k:F_{x_k}\to\bigcup_{y\in D_k}G_y$.  If $\pi\in\S_X$ is any permutation that extends $\pi_k$, then $D_k'$ is a block of $\al\pi\be$.  Again, see Figure \ref{fig:stage1}.

{\bf Stage 4.}  So far, we have defined bijections $\pi_i$ ($i\in I$), $\pi_j$ ($j\in J$), $\pi_k$ ($k\in K$) whose combined (and non-overlapping) domains and codomains are, respectively,
\[
\left( \bigcup_{i\in I} \bigcup_{x\in A_i} E_x \right)
\cup
\left( \bigcup_{j\in J} \bigcup_{x\in C_j} E_x \right)
\cup
\left( \bigcup_{k\in K} F_{x_k} \right)
=
\bigcup_{x\in X}E_x \cup \bigcup_{z\in X_1^1}F_z
\]
and
\[
\left( \bigcup_{i\in I} \bigcup_{y\in B_i} G_y \right)
\cup
\left( \bigcup_{j\in J} H_{x_j} \right)
\cup
\left( \bigcup_{k\in K} \bigcup_{y\in D_k}G_y \right)
=
\bigcup_{y\in X}G_y \cup \bigcup_{z\in X_3^1} H_z.
\]
The complements in $X$ of these sets have cardinality $|X|$, so we may extend $\bigcup_{i\in I}\pi_i \cup \bigcup_{j\in J}\pi_j \cup \bigcup_{k\in K}\pi_k$ arbitrarily to a permutation $\pi\in\S_X$.  By the above discussion, we see that ${\al\pi\be=\ga}$.~\epf

\begin{rem}
The above proof shows that we have the factorization $\P_X=\L_X\S_X\R_X$.  (In fact, since $\S_X\sub\L_X$, we have $\P_X=\L_X\R_X$.)  This is reminiscent of, but quite different to, the factorization of a finite partition monoid as $\P_n=\L_n\I_n\R_n$ utilized in \cite{JEgrpm,JEpnsn}; there, $\I_n$ is (isomorphic to) the symmetric inverse monoid, and the submonoids $\L_n$ and $\R_n$ of $\P_n$ are defined in a very different way to the submonoids $\L_X$ and $\R_X$ of infinite $\P_X$ used here.
\end{rem}


\begin{thm}\label{cor4}
If $X$ is any infinite set, then $\rank(\P_X:\S_X)=2$.
\end{thm}

\pf Proposition \ref{thm3} tells us that $\rank(\P_X:\S_X)\leq2$.  Let $\al\in\P_X$.  The proof will be complete if we can show that $\la\S_X,\al\ra$ is a proper subsemigroup of $\P_X$.  Suppose to the contrary that $\la\S_X,\al\ra=\P_X$.  Let $\be\in\L_X\sm\S_X$, and consider an expression $\be=\ga_1\cdots\ga_r$ where $\ga_1,\ldots,\ga_r\in\S_X\cup\{\al\}$.  Since $\be\not\in\S_X$, at least one of $\ga_1,\ldots,\ga_r$ is equal to $\al$.  Suppose $\ga_1,\ldots,\ga_{s-1}\in\S_X$ but $\ga_s=\al$.  Then $\al\ga_{s+1}\cdots\ga_r=\ga_{s-1}^{-1}\cdots\ga_1^{-1}\be\in\L_X$.  Since $\P_X\sm\L_X$ is a right ideal, it follows that $\al\in\L_X$.  A dual argument shows that $\al\in\R_X$.  But then $\al\in\L_X\cap\R_X=\S_X$, so that $\P_X=\la\S_X,\al\ra=\S_X$, a contradiction. \epf

\begin{rem}
It follows from \cite[Proposition 39]{JEgrpm} and its proof that $\rank(\P_X:\S_X)=2$ for any finite set $X$ with $|X|\geq2$.
\end{rem}

It will be convenient to conclude this section with an extension of Proposition~\ref{thm3}.  The next result shows (among other things) that one of the partitions $\al,\be$ need not have any infinite blocks at all.

\bigskip
\begin{prop}\label{newnewprop}
Let $\al\in\L_X$ and $\be\in\R_X$ be such that $d^*(\al)=d(\be)=|X|$ and either
\bit
	\item[\emph{(i)}] $k^*(\al,2)+d^*(\al,2)=|X|=k(\be,|X|)+d(\be,|X|)$, or
	\item[\emph{(ii)}] $k^*(\al,|X|)+d^*(\al,|X|)=|X|=k(\be,2)+d(\be,2)$.
\eit
Then $\P_X=\la\S_X,\al,\be\ra$.
\end{prop}

\pf Suppose (i) holds.  (The other case is dual.)   By Lemma \ref{lem5file}, we may assume that $d^*(\al,2)=d(\be,|X|)=|X|$.  Write
$
\al=\LXpartn $ and $ \be=\RXpartnt.
$
Note that
\[
\be\al=\partn{C_x}{A_x}{D_x}{B_x}.
\]
We will show that there exists a pair of partitions $\de,\ve\in\la\S_X,\al,\be\ra$ that satisfy the conditions of Proposition \ref{thm3}.
Put
\[
Y=\set{x\in X}{|B_x|\geq2} \AND Z=\set{x\in X}{|D_x|=|X|}.
\]
Let $(Y_1,Y_2)$ and $(Z_1,Z_2)$ be moieties of $Y$ and $Z$, respectively, and let $(U_x)_{x\in X}$ be a moiety of $X$.  Further, suppose $(Y_1^x)_{x\in X}$ is a moiety of $Y_1$.  Write $Z_1=\set{z_x}{x\in X}$ and $Y_1^x=\set{y_w^x}{w\in X}$ for each $x\in X$.
Now let $x\in X$.  For each $w\in X$, choose some $b_w^x\in B_{y_w^x}$.  Choose any embedding $\pi_x^1:A_x\to D_{z_x}$ such that $| D_{z_x}\sm A_x\pi_x^1|=|X|$, and write $D_{z_x}\sm A_x\pi_x^1={\set{d_w^x}{w\in X}}$.  Define $\pi_x^2:\set{b_w^x}{w\in X}\to D_{z_x}\sm A_x\pi_x^1:b_w^x\mt d_w^x$.  Now choose any bijection $\pi_x^3:\bigcup_{w\in X}\left(B_{y_w^x}\sm\{b_w^x\}\right)\to\bigcup_{u\in U_x}C_u$.  Put $\pi_x=\pi_x^1\cup\pi_x^2\cup\pi_x^3$.  So
\[
\pi_x:A_x\cup \bigcup_{y\in Y_1^x}B_y\to D_{z_x}\cup\bigcup_{u\in U_x}C_u
\]
is a bijection.  It is clear that if $\pi\in\S_X$ is any permutation extending $\pi_x$, then $\{x\}\cup \left(\bigcup_{u\in U_x}A_u\right)'$ is a block of $\al\pi\be\al$; see Figure \ref{fig:newnewprop} for an illustration.
\begin{figure}[ht]
   \begin{center}
\begin{tikzpicture}[xscale=1,yscale=1]
  \ublock{4.75-.5}6
  \permblocks{1.5}{5.5}{4.25}{5.5}2
  \permblocks{6}{6.5}{5.5}62
  \blocks{5}{5}{6}{6.5}4
  \lblocks{5}{5.5}4
  \lblocks{4}{4.5}4
  \lblocks{3}{3.5}4
  \permblocks{1.5}{5.5}{-1}{3.75}2
  \draw[dotted,line width=0.4](1.5,4) to [out=270,in=90] (4.25,2);
  \ldotteds{1.5}{2.5+.25}{4.1}
  \block{3-.5}{4.25-.5}{3.5-.5}{4-.5}
  \block{2-.5}{2.5-.5}{2.5-.5}{3-.5}
  \block{1-.5}{1.5-.5}{1.25-.5}{2-.5}
  \ldotted{-.5-.5}{.75-.25}
  \udotteds{-.5-.5}{.5-.5+.25}{-.1}
  \udotteds{4.375}{4.875+.125}{.1}
  \draw(2.125-.5-.25-.125,0.1)node[below]{{\small $\underbrace{\phantom{............................................}}_{{\bigcup_{u\in U_x}A_u}}$}};
  \draw(3.5,4.1)node[above]{{\small $\overbrace{\phantom{.......................................}}^{{\bigcup_{y\in Y_1^x}B_y}}$}};
  \draw(5.125,1.9)node[below]{{\small $\underbrace{\phantom{.................}}_{D_{z_x}}$}};
  \draw(5,6)node[above]{{\small $x$}};
  \draw(6,3)node{{\small $\pi_x^1$}};
  \draw(4.6,2.4)node{{\small $\pi_x^2$}};
  \draw(1.5,2.5)node{{\small $\pi_x^3$}};
  \draw[|-|] (-1-.5,6)--(-1-.5,0);
  \draw[|-|] (-1-.5,4)--(-1-.5,2);
  \draw(-1-.5,5)node[left]{{\small $\al$}};
  \draw(-1-.5,3)node[left]{{\small $\pi$}};
  \draw(-1-.5,1)node[left]{{\small $\be\al$}};
  \permlines{5.375}{5.25+.125}2
  \permlines{5.375-1}{5.25}2
  \permlines{5.375-2}{5.25-.125}2
  \muvs504
	\end{tikzpicture}
    \caption{A schematic diagram of the product $\al\pi\be\al$, focusing on the transversal block $\{x\}\cup\left(\bigcup_{u\in U_x}A_u\right)'$.  See text for further explanation.}
    \label{fig:newnewprop}
   \end{center}
 \end{figure}

The domains of the bijections $\pi_x$ ($x\in X$) are pairwise disjoint, and so too are the codomains.  The complements in $X$ of the domain and codomain of $\bigcup_{x\in X}\pi_x$ have cardinality $|X|$, so we may extend $\bigcup_{x\in X}\pi_x$ arbitrarily to a permutation $\pi\in\S_X$.  By the above discussion, $\{x\}\cup \left(\bigcup_{u\in U_x}A_u\right)'$ is a block of $\ga=\al\pi\be\al$ for all $x\in X$.  It follows that $\ga\in\L_X$, and that $|x|_\ga=|X|$ for all $x\in X$.  By (\ref{lem5.1}.2), we also have $d^*(\ga)\geq d^*(\al)=|X|$.  By Lemma \ref{lem5file}, there exists $\de\in\la\S_X,\ga\ra\cap\L_X$ such that $d^*(\de,|X|)=|X|$ and $|x|_\de=|X|$ for all $x\in X$.  Noting that $k^*(\de,|X|)+d^*(\de,|X|)=|X|=k(\be,2)+d(\be,2)$, a dual argument shows that there exists $\ve\in\la\S_X,\de,\be\ra\sub\la\S_X,\al,\be\ra$ such that $\ve\in\R_X$, $d(\ve,|X|)=|X|$, and $|x'|_\ve=|X|$ for all $x\in X$.  It follows from Proposition \ref{thm3} that $\P_X=\la\S_X,\de,\ve\ra\sub\la\S_X,\al,\be\ra$, and the proof is complete.~\epf


\section{Generating pairs for $\P_X$ modulo $\S_X$}\label{sect:genpairsPXSX}

We saw in Proposition \ref{thm3} that $\P_X$ may be generated by the symmetric group $\S_X$ along with two other partitions.  We call such a pair of partitions a \emph{generating pair} for $\P_X$ modulo $\S_X$.  In this section, we will classify all such generating pairs.  The classification depends crucially on the nature of the cardinal $|X|$, and we will obtain three separate classifications in the cases of $|X|$ being countable (Theorem \ref{thm_countable}), uncountable but regular (Theorem \ref{thm_regular}), and singular (Theorem \ref{thm_singular}).  
We begin with a simple result that will be used in the proof of all three classification theorems.  

\begin{lemma}\label{LRlem}
Suppose $\al,\be\in\P_X$ are such that $\P_X=\la\S_X,\al,\be\ra$.  Then (renaming $\al,\be$ if necessary) $\al\in\L_X$ and $\be\in\R_X$ and $d^*(\al)=d(\be)=|X|$.
\end{lemma}

\pf Consider an expression $\ga=\de_1\cdots\de_r$ where $\ga\in\L_X\sm\S_X$ and $d(\ga)=|X|$, and $\de_1,\ldots,\de_r\in\S_X\cup\{\al,\be\}$.  As in the proof of Theorem \ref{cor4}, it follows that one of $\al,\be$ belongs to $\L_X$.  Without loss of generality, suppose $\al\in\L_X$.  A dual argument shows that one of $\al,\be$ belongs to $\R_X$.  We could not have $\al\in\R_X$ or otherwise $\al\in\S_X$, which would imply that $\P_X=\la\S_X,\be\ra$, contradicting Theorem \ref{cor4}.  So $\be\in\R_X$.  By (\ref{lem5.1}.1), $|X|=d(\ga)\leq d(\de_1)+\cdots+d(\de_r)\leq r\cdot d(\be)$.  It follows that $d(\be)=|X|$.  A dual argument shows that $d^*(\al)=|X|$. \epf

In order to prove our classification theorems, we will need a series of technical lemmas.  The first of these will be used in the proof of all three theorems.


\bigskip
\begin{lemma}\label{prelim-22}
Suppose $\al\in\L_X$ and $\be\in\R_X$, and that $k^*(\al,2)+d^*(\al,2)<|X|$.  If $\ga_1,\ldots,\ga_r\in\S_X\cup\{\al,\be\}$ are such that $\ga_1\cdots\ga_r\in\L_X$, then $k^*(\ga_1\cdots\ga_r,2)+d^*(\ga_1\cdots\ga_r,2)<|X|$.
\end{lemma}

\pf Write $\al=\partnL{A_x}{B_i}$ and $\be=\partnR{C_x}{D_j}$.  
The result is clearly true if $r=1$, so suppose $r\geq2$, and put $\ga=\ga_1\cdots\ga_{r-1}$.  Since $\P_X\sm\L_X$ is a right ideal, it follows that $\ga\in\L_X$.  Thus, an induction hypothesis gives $k^*(\ga,2)+d^*(\ga,2)<|X|$.  
Write $\ga=\partnL{E_x}{F_k}$.
We now break the proof up into three cases.

{\bf Case 1.}  If $\ga_r\in\S_X$, then clearly $k^*(\ga\ga_r,2)=k^*(\ga,2)$ and $d^*(\ga\ga_r,2)=d^*(\ga,2)$, so the inductive step is trivial in this case.

{\bf Case 2.}  Next suppose $\ga_r=\al$.  Since $2$ is a regular cardinal, $k^*(\ga\al,2)\leq k^*(\ga,2)+k^*(\al,2)<|X|$ by (\ref{k-lemma}.2).  We also have $d^*(\ga\al,2)\leq d^*(\ga,2)+d^*(\al,2)+k^*(\al,2)<|X|$ by (\ref{d-lemma}.2).




{\bf Case 3.}  Finally, suppose $\ga_r=\be$.
Write $\ga\be=\partnL{P_x}{Q_l}$.
Let
\[
Y=\set{x\in X}{|P_x|\geq2} \AND M=\set{l\in L}{|Q_l|\geq2}.
\]
We must show that $|Y|<|X|$ and $|M|<|X|$.  We begin with $Y$.  Put
$
Y_1 = \set{x\in Y}{|E_x|=1} $ and $ Y_2 = \set{x\in Y}{|E_x|\geq2}.
$
Now $|Y_2|\leq k^*(\ga,2)<|X|$ so, to show that $|Y|<|X|$, it remains to show that $|Y_1|<|X|$.  Now suppose $x\in Y_1$, and write $E_x=\{e_x\}$.  We claim that there exists $k_x\in K$ such that $|F_{k_x}|\geq2$ and $(u_x,e_x)\in\ker(\be)$ for some $u_x\in F_{k_x}$.  The proof of the claim breaks up into two subcases.

{\bf Subcase 3.1.}  First suppose that $e_x\in C_y$ for some $y\in X$.  If $|C_y|=1$, then we would have $P_x=\{y\}$, contradicting the fact that $|P_x|\geq2$.  So $|C_y|\geq2$.  Now $C_y\sm\{e_x\}$ has trivial intersection with $E_z$ for each $z\in X\sm\{x\}$, or else then we would have $(x,z)\in\ker(\ga\be)$ for some $z\not=x$, contradicting the fact that $\ga\be\in\L_X$.  So, for all $u\in C_y\sm\{e_x\}$, we have $u\in F_{k_u}$ for some $k_u\in K$.  (The map $u\mt k_u$ need not be injective.)  If $|F_{k_u}|=1$ for all $u\in C_y\sm\{e_x\}$ then, again, we would have $P_x=\{y\}$, a contradiction.  So it follows that $|F_{k_u}|\geq2$ for at least one $u\in C_y\sm\{e_x\}$.  We now choose $u_x$ to be any such $u$, and we put $k_x=k_{u_x}$.

{\bf Subcase 3.2.}  The case in which $e_x\in D_y$ for some $y\in X$ is similar to the previous subcase.

With the claim established, we note that the map $Y_1\to K:x\mt k_x$ is injective.  Indeed, if $k_{x_1}=k_{x_2}$ for some $x_1,x_2\in Y_1$, then we would have $(x_1,x_2)\in\ker(\ga\be)$, which implies that $x_1=x_2$.  Now the image of $Y_1$ under this map is contained in the set $\set{k\in K}{|F_k|\geq2}$ which has cardinality $d^*(\ga,2)$.  Thus, $|Y_1|\leq d^*(\ga,2)<|X|$ as required.

Next suppose $l\in M$.  Fix some $x\in Q_l$.  Now $C_x$ has trivial intersection with $E_y$ for each $y\in X$ (or else $Q_l'$ would not be a nontransversal block of $\ga\be$).  Let $N=\set{k\in K}{F_k\cap C_x\not=\emptyset}$.  So $N\not=\emptyset$.  If $|F_k|=1$ for all $k\in N$, then we would have $Q_l=\{x\}$, contradicting the fact that $|Q_l|\geq2$.  So there exists some $n_l\in N$ such that $|F_{n_l}|\geq2$.  Again, the map $M\to\set{k\in K}{|F_k|\geq2}:l\mt n_l$ is injective, so it follows that $|M|\leq d^*(\ga,2)<|X|$, as required.  This completes the inductive step in Case 3. \epf

The next lemma will be of use in the case that $X$ is uncountable, whether regular or singular.

\bigskip
\begin{lemma}\label{prelim-nu}
Suppose $\ao\leq\mu\leq|X|$ is a regular cardinal.  Suppose $\al\in\L_X$ and $\be\in\R_X$ are such that $k^*(\al,\mu)+d^*(\al,\mu)<|X|$ and $k(\be,\mu)+d(\be,\mu)<|X|$.  If $\ga_1,\ldots,\ga_r\in\S_X\cup\{\al,\be\}$ are such that $\ga_1\cdots\ga_r\in\L_X$, then ${k^*(\ga_1\cdots\ga_r,\mu)+d^*(\ga_1\cdots\ga_r,\mu)<|X|}$.
\end{lemma}

\pf Write $\al=\partnL{A_x}{B_i}$ and $\be=\partnR{C_x}{D_j}$.
Again, the $r=1$ case is trivial, so suppose $r\geq2$, and put $\ga=\ga_1\cdots\ga_{r-1}\in\L_X$.  An induction hypothesis gives $k^*(\ga,\mu)+d^*(\ga,\mu)<|X|$.  Write $\ga=\partnL{E_x}{F_k}$.
We now break the proof up into three cases.

{\bf Case 1.}  The $\ga_r\in\S_X$ case is trivial.

{\bf Case 2.}  Again, the case in which $\ga_r=\al$ follows from (\ref{k-lemma}.2) and (\ref{d-lemma}.2).




{\bf Case 3.}  Finally, suppose $\ga_r=\be$.  Write $\ga\be=\partnL{P_x}{Q_l}$.
Let
\[
Y=\set{x\in X}{|P_x|\geq\mu} \AND M=\set{l\in L}{|Q_l|\geq\mu}.
\]
We must show that $|Y|<|X|$ and $|M|<|X|$.  We begin with~$Y$.  Put
$
Y_1 = \set{x\in Y}{|E_x|<\mu} $ and $ Y_2 = \set{x\in Y}{|E_x|\geq\mu}.
$
Now $|Y_2|\leq k^*(\ga,\mu)<|X|$ so, to show that $|Y|<|X|$, it remains to show that $|Y_1|<|X|$.  Now suppose $x\in Y_1$.  Consider the connected component containing $x$ in the product graph $\Ga(\ga,\be)$.  The middle row of this connected component is (omitting double dashes, for convenience)
\[
E_x \cup \bigcup_{k\in K_x}F_k = \bigcup_{y\in P_x}C_y \cup \bigcup_{z\in Z_x}D_z
\]
for some subsets $K_x\sub K$ and $Z_x\sub X$.  We now claim that one of the following holds: (i)~$|F_k|\geq\mu$ for some $k\in K_x$, (ii) $|C_y|\geq\mu$ for some $y\in P_x$, or (iii) $|D_z|\geq\mu$ for some $z\in Z_x$.  Indeed, suppose not.  Put
\[
\G = \{E_x\} \cup \set{F_k}{k\in K_x} \AND \H = \set{C_y}{y\in P_x}\cup\set{D_z}{z\in Z_x}.
\]
Note that $|G|<\mu$ and $|H|<\mu$ for all $G\in\G$ and $H\in\H$.  Fix some $a\in P_x$.  For any $b\in P_x$, there exists a sequence of points $c_1,c_2,\ldots,c_{2s}\in X$ such that 
\begin{align*}
c_1 &\in C_a\cap G_1 &&\text{for some $G_1\in\G$}\\
c_2 &\in G_1\cap H_1 &&\text{for some $H_1\in\H$}\\
c_3 &\in H_1\cap G_2 &&\text{for some $G_2\in\G$}\\
c_4 &\in G_2\cap H_2 &&\text{for some $H_2\in\H$}\\
& \ \:\!\:\!\:\! \vdots \\
c_{2s-2} &\in G_{s-1}\cap H_{s-1} &&\text{for some $H_{s-1}\in\H$}\\
c_{2s-1} &\in H_{s-1}\cap G_s &&\text{for some $G_s\in\G$}\\
c_{2s} &\in G_s\cap C_b.
\end{align*}
Let $\xi$ denote the number of such sequences, and let $\xi_s$ be the number of such sequences for fixed $s$.  Let us give an upper bound for $\xi_s$.  There are at most $|C_a|$ choices for $c_1$.  Once $c_1$ is chosen, $G_1$ is uniquely determined, and then there are at most $|G_1|$ choices for $c_2$.  Once $c_2$ is chosen, $H_1$ is uniquely determined, and then there are at most $|H_1|$ choices for $c_3$.  Continuing in this fashion, we see that $\xi_s\leq|C_a|\times|G_1|\times|H_1|\times\cdots\times|H_{s-1}|\times|G_s|<\mu^{2s}=\mu$.  It follows that $\xi=\xi_1+\xi_2+\xi_3+\cdots<\mu$ since $\mu\geq\aleph_1$ is regular.  But this implies that there are less than $\mu$ choices for $b\in P_x$.  That is, $|P_x|<\mu$, a contradiction.  So, indeed, one of (i), (ii), (iii) must hold.  But this implies that $|Y_1|\leq d^*(\ga,\mu)+k(\be,\mu)+d(\be,\mu)<|X|$.  This completes the proof that $|W|<|X|$.

A similar argument shows that $|M|\leq d^*(\ga,\mu)+k(\be,\mu)+d(\be,\mu)<|X|$.  This completes the inductive step in Case 3. \epf

\begin{rem}
The argument used in the last stage of Case 3 in the above proof is reminiscent of the proof of \cite[Lemma 27]{JEifims}.  Everything in the above proof works for $\mu=\an$ except for the final claim that $\xi_1+\xi_2+\xi_3+\cdots<\mu=\an$, which is not the case if infinitely many of the $\xi_s$ are nonzero.  It is this fact that will enable us to generate $\P_X$, for countable $X$, with $\S_X$ along with two additional partitions, neither of which have any infinite blocks; see Proposition \ref{newnewnewprop} below.
\end{rem}

The next result shows that the converse of Proposition \ref{newnewprop} holds in the case of $|X|$ being regular but uncountable.

\bigskip
\begin{thm}\label{thm_regular}
Suppose $|X|\geq\aleph_1$ is a regular cardinal.  Then $\P_X=\la\S_X,\al,\be\ra$ if and only if (up to renaming $\al,\be$ if necessary) $\al\in\L_X$, $\be\in\R_X$, $d^*(\al)=d(\be)=|X|$ and either
\bit
	\item[\emph{(i)}] $k^*(\al,2)+d^*(\al,2)=|X|=k(\be,|X|)+d(\be,|X|)$, or
	\item[\emph{(ii)}] $k^*(\al,|X|)+d^*(\al,|X|)=|X|=k(\be,2)+d(\be,2)$.
\eit
\end{thm}

\pf The reverse implication was proved in Proposition \ref{newnewprop}.  So suppose $\P_X=\la\S_X,\al,\be\ra$.  By Lemma \ref{LRlem}, and renaming $\al,\be$ if necessary, we may assume that $\al\in\L_X$, $\be\in\R_X$, and $d^*(\al)=d(\be)=|X|$.  
Suppose that (i) and (ii) do not hold.  So one of
\bit
	\item[(I)] $k^*(\al,2)+d^*(\al,2)<|X|$, or
	\item[(II)] $k(\be,|X|)+d(\be,|X|)<|X|$
\eit
holds, and so too does one of
\bit
	\item[(III)] $k^*(\al,|X|)+d^*(\al,|X|)<|X|$, or
	\item[(IV)] $k(\be,2)+d(\be,2)<|X|$.
\eit
We will show that we obtain a contradiction in any case.  If (I) holds, then Lemma \ref{prelim-22} tells us that $\la\S_X,\al,\be\ra$ does not contain any $\ga\in\L_X$ with $k^*(\ga,2)=|X|$, contradicting the assumption that $\P_X=\la\S_X,\al,\be\ra$.  Dually, (IV) leads to a contradiction too.  If (II) and (III) both hold, then Lemma \ref{prelim-nu} (with $\mu=|X|$) tells us that $\la\S_X,\al,\be\ra$ does not contain any $\ga\in\L_X$ with $k^*(\ga,|X|)=|X|$, a contradiction.  \epf

Now that we have achieved a classification in the case of $X$ being regular but uncountable, we move on to the countable case.  The next result shows that if $X$ is countable, then $\P_X$ may be generated by $\S_X$ along with two partitions, neither of which has any infinite blocks.

\bigskip
\begin{prop}\label{newnewnewprop}
Suppose $|X|=\an$.  Let $\al\in\L_X$ and $\be\in\R_X$ be such that $d^*(\al)=d(\be)=\an$ and either
\bit
	\item[\emph{(i)}] $k^*(\al,2)+d^*(\al,2)=\an=k(\be,2)+d(\be,3)$, or
	\item[\emph{(ii)}] $k^*(\al,2)+d^*(\al,3)=\an=k(\be,2)+d(\be,2)$.
\eit
Then $\P_X=\la\S_X,\al,\be\ra$.
\end{prop}

\pf Suppose (ii) holds.  (The other case is dual.)  By Proposition \ref{newnewprop}, it is enough to show that there exists $\ga\in\la\S_X,\al,\be\ra$ such that $\ga\in\L_X$ and $k^*(\ga,\an)=d^*(\ga)=\an$.  By the dual of Lemma \ref{lem5file}, we may assume that $d(\be,2)=\an$.  Write
$
\al=\LXpartn $ and $ \be=\RXpartnt.
$
We first claim that there exists $\de\in\la\S_X,\al\ra$ such that $\de\in\L_X$ and $d^*(\de,3)=\an$.  Indeed, we put $\de=\al$ if $d^*(\al,3)=\an$.  Otherwise, suppose $k^*(\al,2)=\an$.  Let $(Y_1,Y_2)$ be a moiety of $Y=\set{x\in X}{|A_x|\geq2}$.  Since $|X\sm Y_1|=\an=\left|X\sm\bigcup_{x\in Y_1}A_x\right|$, we may choose a permutation $\pi\in\S_X$ such that $\left(\bigcup_{x\in Y_1}A_x\right)\pi=Y_1$.  Then $\al\pi\al\in\L_X$ and, for each $x\in Y_1$, $\{x\}\cup\left(\bigcup_{y\in A_x}A_{y\pi}\right)'$ is a block of $\al\pi\al$, and $\left|\bigcup_{y\in A_x}A_{y\pi}\right|\geq4$.  That is, $k^*(\al\pi\al,4)=\an$.  But then, by Lemma \ref{lem5file}, there exists $\de\in\la\S_X,\al\pi\al\ra\sub\L_X$ with $d^*(\de,4)\geq k^*(\al\pi\al,4)=\an$.  This completes the proof of the claim, since $d^*(\de,3)\geq d^*(\de,4)=\an$.  Now write 
$\de=\LXpartnth$.
Let $U=\set{x\in X}{|F_x|\geq3}$ and $V=\set{x\in X}{|D_x|\geq2}$.  Note that
\[
\be\de = \partn{C_x}{E_x}{D_x}{F_x}.
\]
Our goal will be to construct a permutation $\si\in\S_X$ such that $\ga=\de\si\be\de\in\L_X$ and ${k^*(\ga,\an)=\an}$.  The proof will then be complete since we will also have $d^*(\ga)\geq d^*(\de)=\an$.  
Now let $\infty$ be a symbol that is not an element of $X$.
Let $(U_x)_{x\in X\cup\{\infty\}}$, $(V_x)_{x\in X\cup\{\infty\}}$ and $(Z_x)_{x\in X}$ be moieties of $U$, $V$ and $X$, respectively.
For each $x\in X$, write $U_x=\set{u_x^r}{r\in\N}$ and $V_x=\set{v_x^r}{r\in\N}$.  Now fix $x\in X$.  We define a bijection
\[
\si_x:E_x\cup\bigcup_{y\in U_x}F_y \to \bigcup_{y\in V_x}D_y \cup \bigcup_{z\in Z_x}C_z
\]
as follows.  First, choose some $e_x\in E_x$ and $a_x^r,b_x^r\in F_{u_x^r}$, $c_x^r,d_x^r\in D_{v_x^r}$ where $a_x^r\not=b_x^r$ and $c_x^r\not=d_x^r$ for each $r\in \N$.  We then define a bijection 
\[
\si_x^1:\{e_x\}\cup\set{a_x^r,b_x^r}{r\in \N}\to\set{c_x^r,d_x^r}{r\in \N}
\]
by $e_x\si_x^1=c_x^1$, $a_x^r\si_x^1=d_x^r$ and $b_x^r\si_x^1=c_x^{r+1}$ for each $r\in \N$.  Since $|F_y|\geq3$ for all $y\in U_x$ and since $|Z_x|=\an$, we see that the complements of the domain and codomain of $\si_x^1$
in $E_x\cup\bigcup_{y\in U_x}F_y$ and $\bigcup_{y\in V_x}D_y \cup \bigcup_{z\in Z_x}C_z$ (respectively) have cardinality~$\an$.  So we extend $\si_x^1$ arbitrarily to a bijection $\si_x:E_x\cup\bigcup_{y\in U_x}F_y\to\bigcup_{y\in V_x}D_y \cup \bigcup_{z\in Z_x}C_z$.  It is clear that if $\si\in\S_X$ is any permutation that extends $\si_x$, then $\{x\}\cup\left(\bigcup_{z\in Z_x}E_z\right)'$ is a block of $\de\si\be\de$; see Figure \ref{fig:newnewnewprop}.
\begin{figure}[ht]
   \begin{center}
\begin{tikzpicture}[xscale=1,yscale=1]
  \blocks{5}{5}{6}{6.5}4
  \lblocks{5}{5.5}4
  \lblocks{4}{4.5}4
  \lblocks{3}{3.5}4
  \ublock{5-.5}{5.5-.5}
  \ublock{4-.5}{4.5-.5}
  \ublock{3-.5}{3.5-.5}
  \permblocks{1.5}{6.5}{1}{10}2
  \block{5.5}{6.5}{5.5}{6}
  \block{7}{7.5}{6.5}{6.5}
  \block{8}{8.5}{7}{8}
  \udotteds{8.75}{10}{-.1}
  \ldotted{8.25}{9.5}
  \udotteds{1.5-.5}{2.5+.25-.5}{-.1}
  \udotteds{1.5}{2.75}{2.1}
  \draw(5.25,4)node[above]{{\small $\overbrace{\phantom{....}}^{{F_{u_x^1}}}$}};
  \draw(4.25,4)node[above]{{\small $\overbrace{\phantom{....}}^{{F_{u_x^2}}}$}};
  \draw(3.25,4)node[above]{{\small $\overbrace{\phantom{....}}^{{F_{u_x^3}}}$}};
  \draw(4.75,2)node[below]{{\small $\underbrace{\phantom{....}}_{{D_{v_x^1}}}$}};
  \draw(3.75,2)node[below]{{\small $\underbrace{\phantom{....}}_{{D_{v_x^2}}}$}};
  \draw(2.75,2)node[below]{{\small $\underbrace{\phantom{....}}_{{D_{v_x^3}}}$}};
  \draw(7.5,0.1)node[below]{{\small $\underbrace{\phantom{.......................................}}_{{\bigcup_{z\in Z_x}E_z}}$}};
  \draw(5,6)node[above]{{\small $x$}};
  \draw[|-|] (0,6)--(0,0);
  \draw[|-|] (0,4)--(0,2);
  \draw(0,5)node[left]{{\small $\de$}};
  \draw(0,3)node[left]{{\small $\si$}};
  \draw(0,1)node[left]{{\small $\be\de$}};
  \permlines{5.125+1}{4.875}2
  \permlines{5.375}{4.75}2
  \permlines{5.375-.125}{3.875}2
  \permlines{5.375-1}{4.75-1}2
  \permlines{5.375-.125-1}{3.875-1}2
  \permlines{5.375-2}{4.75-2}2
  \muvs{3.25}02 \draw (3.25,4) to [out=270,in=45] (3.25-.7125,3);
  \muvs504
  \draw(2,3)node{{\small $\si_x$}};
	\end{tikzpicture}
    \caption{A schematic diagram of the product $\de\si\be\de$, focusing on the transversal block $\{x\}\cup\left(\bigcup_{z\in Z_x}E_z\right)'$.  See text for further explanation.}
    \label{fig:newnewnewprop}
   \end{center}
 \end{figure}

Now, $\bigcup_{x\in X}\si_x$ has domain 
\[
\bigcup_{x\in X}E_x\cup\bigcup_{x\in X}\bigcup_{y\in U_x}F_y = 
\bigcup_{x\in X}E_x\cup\bigcup_{y\in U\sm U_\infty}F_y
\]
and codomain
\[
\bigcup_{x\in X}\bigcup_{y\in V_x}D_y\cup\bigcup_{x\in X}\bigcup_{z\in Z_x}C_z = 
\bigcup_{y\in V\sm V_\infty}D_y\cup\bigcup_{x\in X}C_x.
\]
The complements in $X$ of these sets 
have cardinality $\an$.  So we extend $\bigcup_{x\in X}\si_x$ arbitrarily to $\si\in\S_X$.  By the above discussion, $\{x\}\cup\left(\bigcup_{z\in Z_x}E_z\right)'$ is a block of $\ga=\de\si\be\de$ for each $x\in X$.  It follows that $\ga\in\L_X$.  Also, since $|Z_x|=\an$, it follows that $|x|_\ga=\an$ for all $x\in X$, and so $k^*(\ga,\an)=\an$, as required.  This completes the proof. \epf

We need the next lemma to show that the converse of Proposition \ref{newnewnewprop} is true.

\bigskip
\begin{lemma}\label{prelim-23}
Suppose $|X|=\an$ and $\al\in\L_X$ and $\be\in\R_X$ are such that $k^*(\al,2)+d^*(\al,3)<\an$ and $k(\be,2)+d(\be,3)<\an$.  If $\ga_1,\ldots,\ga_r\in\S_X\cup\{\al,\be\}$ are such that $\ga_1\cdots\ga_r\in\L_X$, then $k^*(\ga_1\cdots\ga_r,2)+d^*(\ga_1\cdots\ga_r,3)<\an$.
\end{lemma}

\pf Write $\al=\partnL{A_x}{B_i}$ and $\be=\partnR{C_x}{D_j}$.  Again, the $r=1$ case is trivial, so suppose $r\geq2$, and put $\ga=\ga_1\cdots\ga_{r-1}\in\L_X$.  An induction hypothesis gives $k^*(\ga,2)+d^*(\ga,3)<\an$.  Write $\ga=\partnL{E_x}{F_k}$.
We now break the proof up into three cases.

{\bf Case 1.}  The $\ga_r\in\S_X$ case is trivial.

{\bf Case 2.}  Next suppose $\ga_r=\al$.  Again, (\ref{k-lemma}.2) gives $k^*(\ga\al,2)\leq k^*(\ga,2)+k^*(\al,2)<\an$.
But~$3$ is not a regular cardinal, so (\ref{d-lemma}.2) is not of any use here.  Now $\ga\al = \partnL{G_x}{B_i,H_k}$, where
$G_x=\bigcup_{y\in E_x}A_y$ and $H_k = \bigcup_{y\in F_k}A_y$
for each $x\in X$ and $k\in K$.
%
%
%
There are $d^*(\al,3)$ values of $i\in I$ such that $|B_i|\geq3$.  Next suppose $k\in K$ is such that $|H_k|\geq3$.  Then either (i)~$|F_k|\geq3$, (ii) $|F_k|=2$ and $|A_y|\geq2$ for some $y\in F_k$, or (iii) $|F_k|=1$ and $|A_y|\geq3$ where $F_k=\{y\}$.  There are $d^*(\ga,3)$ values of $k$ for which (i) holds, at most $k^*(\al,2)$ values of $k$ for which (ii) holds, and at most $k^*(\al,3)\leq k^*(\al,2)$ values of $k$ for which (iii) holds.  Thus, $d^*(\ga\al,3)\leq d^*(\al,3)+d^*(\ga,3)+2k^*(\al,2)<\an$.  This completes the inductive step in this case.

{\bf Case 3.}  Finally, suppose $\ga_r=\be$.  Write $\ga\be=\partnL{P_x}{Q_l}$.
Let
\[
Y=\set{x\in X}{|P_x|\geq2} \AND M=\set{l\in L}{|Q_l|\geq3}.
\]
We must show that $|Y|<\an$ and $|M|<\an$.  Put
$
Y_1 = \set{x\in Y}{|E_x|=1} $ and $ Y_2 = \set{x\in Y}{|E_x|\geq2}.
$
Now $|Y_2|\leq k^*(\ga,2)<\an$ so, to show that $|Y|<\an$, it remains to show that $|Y_1|<\an$.  Now suppose $x\in Y_1$.  Consider the connected component containing $x$ in the product graph $\Ga(\ga,\be)$.  The middle row of this connected component is (omitting double dashes)
\[
E_x \cup \bigcup_{k\in K_x}F_k = \bigcup_{y\in P_x}C_y \cup \bigcup_{z\in Z_x}D_z
\]
for some subsets $K_x\sub K$ and $Z_x\sub X$.  We now claim that one of the following holds: (i)~$|F_k|\geq3$ for some $k\in K_x$, (ii) $|C_y|\geq2$ for some $y\in P_x$, or (iii) $|D_z|\geq3$ for some $z\in Z_x$.  Indeed, suppose not.  Write $C_y=\{c_y\}$ for each $y\in P_x$.  Choose some $a\in P_x$.  Now, if $c_a\in E_x$, then $\{x,a'\}$ would be a block of $\ga\be$ (since $|E_x|=1$), contradicting the fact that $|P_x|\geq2$.  So $c_a\in F_{k_1}$ for some $k_1\in K_x$.  If $|F_{k_1}|=1$, then $\{a'\}$ would be a block of $\ga\be$, a contradiction.  So $F_{k_1}=\{c_a,w_1\}$ for some $w_1\in X\sm\{c_a\}$.  If $w_1\in C_b$ for some $b\in P_x$, then $\{a',b'\}$ would be a block of $\ga\be$, a contradiction.  So we must have $w_1\in D_{z_1}$ for some $z_1\in Z_x$.  
If $|D_{z_1}|=1$, then $\{a'\}$ would be a block of $\ga\be$, so we must have $D_{z_1}=\{w_1,w_2\}$ for some $w_2\in X\sm\{c_a,w_1\}$.  Continuing in this way, we see that the connected component of $a'$ in the product graph $\Ga(\ga,\be)$ is $\{a',c_a'',w_1'',w_2'',w_3'',\ldots\}$ for some $w_1,w_2,w_3,\ldots\in X$.  But this says that $\{a'\}$ is a block of $\ga\be$, contradicting the fact that $x\in[a']_{\ga\be}$.  
So, indeed, one of (i), (ii), (iii) must hold.  But this implies that $|Y_1|\leq d^*(\ga,3)+k(\be,2)+d(\be,3)<\an$.


A similar argument shows that $|M|\leq d^*(\ga,3)+k(\be,2)+d(\be,3)<\an$.  This completes the inductive step in Case 3. \epf

The proof of the following is virtually identical to the proof of Theorem \ref{thm_regular}.  Instead of applying Lemmas \ref{prelim-22} and \ref{prelim-nu}, we apply Lemmas \ref{prelim-22} and \ref{prelim-23}.

\bigskip
\begin{thm}\label{thm_countable}
Suppose $|X|=\an$ and let $\al,\be\in\P_X$.  Then $\P_X=\la\S_X,\al,\be\ra$ if and only if (renaming $\al,\be$ if necessary), $\al\in\L_X$, $\be\in\R_X$, $d^*(\al)=d(\be)=\an$ and either
\bit
	\item[\emph{(i)}] $k^*(\al,2)+d^*(\al,2)=\an=k(\be,2)+d(\be,3)$, or
	\item[\emph{(ii)}] $k^*(\al,2)+d^*(\al,3)=\an=k(\be,2)+d(\be,2)$. \epfres
\eit
\end{thm}


We now turn our attention to the case of $X$ having singular cardinality.

\bigskip
\begin{prop}\label{newnewnewnewprop}
Suppose $|X|$ is singular.  Let $\al\in\L_X$ and $\be\in\R_X$ be such that $d^*(\al)=d(\be)=|X|$ and either
\bit
	\item[\emph{(i)}] $k^*(\al,2)+d^*(\al,2)=|X|=k(\be,\mu)+d(\be,\mu)$ for all cardinals $\mu<|X|$, or
	\item[\emph{(ii)}] $k^*(\al,\mu)+d^*(\al,\mu)=|X|=k(\be,2)+d(\be,2)$ for all cardinals $\mu<|X|$.
\eit
Then $\P_X=\la\S_X,\al,\be\ra$.
\end{prop}

\pf Suppose (ii) holds.  (The other case is dual.)  By Proposition \ref{newnewprop}, it is enough to show that there exists $\ga\in\la\S_X,\al,\be\ra$ such that $\ga\in\L_X$ and $k^*(\ga,|X|)=d^*(\ga)=|X|$.  Since $|X|$ is singular, we have $X=\bigcup_{i\in I}X_i$, where $|I|<|X|$ and $|X_i|<|X|$ for all $i\in I$.  Put $\ka=|I|$ and $\lam_i=|X_i|$ for each $i\in I$.  Write
$\al=\LXpartn$ and $\be=\RXpartnt$.
In order to avoid notational ambiguity, we will suppose that $\infty$ is a symbol that does not belong to $X$, and that $I$ does not contain $1$ or $2$.

By the dual of Lemma \ref{lem5file}, we may suppose without loss of generality that $d(\be,2)=|X|$.  We claim that either
\bit
\item[(a)] $k^*(\al,\mu)=|X|$ for all cardinals $\mu<|X|$, or
\item[(b)] $d^*(\al,\mu)=|X|$ for all cardinals $\mu<|X|$.
\eit
Indeed, suppose (a) does not hold.  Then there exists $\nu<|X|$ such that $k^*(\al,\nu)<|X|$.  But then, for any $\nu\leq\mu<|X|$, $k^*(\al,\mu)\leq k^*(\al,\nu)<|X|$ which, together with $k^*(\al,\mu)+d^*(\al,\mu)=|X|$, implies that $d^*(\al,\mu)=|X|$ holds for all $\nu\leq\mu<|X|$.  But also, for any $\mu<\nu$, $d^*(\al,\mu)\geq d^*(\al,\nu)=|X|$, so it follows that (b) holds.  We will now break the proof up into two cases, according to whether (a) or (b) holds.

{\bf Case 1.}  Suppose (a) holds.  For each cardinal $\mu\leq|X|$, let $Y_\mu=\set{x\in X}{|A_x|\geq\mu}$.  Note that $|Y_\mu|=|X|$ for all $\mu<|X|$, and that $Y_\mu\sub Y_\nu$ if $\nu\leq\mu\leq|X|$.  For each $x\in Y_\ka$, choose a subset $E_x\sub A_x$ with $|E_x|=\ka$ and $A_x\sm E_x\not=\emptyset$,
and write $E_x=\set{e_{xi}}{i\in I}$.

Next, suppose $(X_1,X_2)$ is a moiety of $X$.  We claim that either (1) $|X_1\cap Y_\mu|=|X|$ for all $\mu<|X|$, or (2) $|X_2\cap Y_\mu|=|X|$ for all $\mu<|X|$.  Indeed, suppose (2) does not hold.  Then there exists $\nu<|X|$ such that $|X_2\cap Y_\nu|<|X|$.  Then, for any $\nu\leq\mu<|X|$, $|X_2\cap Y_\mu|\leq|X_2\cap Y_\nu|<|X|$.  But for any $\mu\leq|X|$, $Y_\mu=(X_1\cap Y_\mu)\sqcup(X_2\cap Y_\mu)$, so it follows that $|X_1\cap Y_\mu|=|X|$ for all $\nu\leq\mu<|X|$.  If $\mu<\nu$, then $|X_1\cap Y_\mu|\geq|X_1\cap Y_\nu|=|X|$.  So (1) holds, and the claim is established.  From now on, we fix $X_1,X_2$ as above and suppose, without loss of generality, that (1) holds.

We will now construct, by transfinite recursion, a set $Z=\set{d_{xi}}{x\in Y_\ka,i\in I}\sub X_1$ such that $d_{xi}\in Y_{\lam_i}$ for each $(x,i)\in Y_\ka\times I$.  Indeed, fix some well-ordering $<$ on $Y_\ka\times I$.  Suppose $(x,i)\in Y_\ka\times I$ and that we have already defined the elements $Z_{xi}=\set{d_{yj}}{(y,j)<(x,i)}$.  Since this set is constructed recursively, by adding a single element at a time, we see that $|Y_{\lam_i}\sm Z_{xi}|=|Y_{\lam_i}|=|X|$.  So we define $d_{xi}$ to be any element of $Y_{\lam_i}\sm Z_{xi}$.
With $Z$ so defined, there is a natural bijection $\si:\bigcup_{x\in Y_\ka} E_x\to Z:e_{xi}\mt d_{xi}$.  Since the complement in $X$ of the domain and codomain of $\si$ both have cardinality $|X|$, we may extend $\si$ arbitrarily to a permutation $\pi\in\S_X$.  Now put $\ga=\al\pi\al\in\L_X$.
Then $\ga=\partnL{F_x}{B_x,G_x}$,
where $F_x=\bigcup_{y\in A_x}A_{y\pi}$ and $G_x=\bigcup_{y\in B_x}A_{y\pi}$ for each $x\in X$.  Let $x\in Y_\ka$.  Then $A_{d_{xi}}=A_{e_{xi}\pi}\sub F_x$ for all $i\in I$.  So $F_x\sups \bigcup_{i\in I}A_{d_{xi}}$, and $|F_x|\geq\sum_{i\in I}|A_{d_{xi}}|\geq\sum_{i\in I}\lam_i=|X|$.  Since $|Y_\ka|=|X|$, it follows that $k^*(\ga,|X|)=|X|$, as required. 

{\bf Case 2.}  Now suppose (b) holds.  By the previous case, it is sufficient to show that there exists $\de\in\la\S_X,\al,\be\ra$ such that $\de\in\L_X$, $d^*(\de)=|X|$ and $k^*(\de,\mu)=|X|$ for all cardinals $\mu<|X|$.

This time, for each cardinal $\mu\leq|X|$, we define $Z_\mu=\set{x\in X}{|B_x|\geq\mu}$.  Let $(W_1,W_2)$ be a moiety of $Z_{\an}=\set{x\in X}{|B_x|\geq\an}$.  As in the previous case, we may assume that $|W_1\cap Z_\mu|=|X|$ for all $\mu<|X|$.  Write $W_1=\set{w_x}{x\in X}$, and let $(W_2^x)_{x\in X\cup\{\infty\}}$ be a moiety of $W_2$.  Let $(U_1,U_2)$ be a moiety of $\set{x\in X}{|D_x|\geq2}$ and write $U_1=\set{u_x}{x\in X}$.  For each $x\in X$, choose $a_x\in A_x$, $b_x,c_x\in D_{u_x}$, $d_x\in B_{w_x}$ with $b_x\not=c_x$.  Let $(V_x)_{x\in X}$ be a moiety of $X$ and, for each $x\in X$, let $V_x=V_x^1\sqcup V_x^2$ where $|V_x^1|=|B_{w_x}|$ and $|V_x^2|=|X|$.  Write $V_x^1=\set{v_{xy}}{y\in B_{w_x}\sm\{d_x\}}$, noting that $|B_{w_x}|\geq\an$.  For each $y\in B_{w_x}\sm\{d_x\}$, choose some $e_{xy}\in C_{v_{xy}}$.

\begin{figure}[ht]
   \begin{center}
\begin{tikzpicture}[xscale=1,yscale=1]
  \blocks{5}{5}{6}{6.5}4
  \lblocks{3.5}{5.5}4
  \lblocks7{7.5}4
  \lblocks8{8.5}4
  \lblocks9{9.5}4
  \block{6.5}7{6.5}{7.5}
  \block{7.5}88{8.5}
  \block{8.5}{9}{9}{10}
  \ublock{5}6
  \permblocks{3.5}{11}{-.25}{10.5}2
  \block{3-.5+.75}{4.25-.5+.75}{3.5-.5+.75}{4-.5+.75}
  \block{2-.5+.75}{2.5-.5+.75}{2.5-.5+.75}{3-.5+.75}
  \block{1-.5+.75}{1.5-.5+.75}{1.25-.5+.75}{2-.5+.75}
  \permlines{5.375+.75}{5.25+.125+.5}2
  \permlines{5.375}{5.125}2
  \permlines{5.375-.125}{4.375}2
  \permlines{5.375-.125-.125}{2.625}2
  \lpermlinels{5.375-.125-.125-.125}{1.625}2
  \ldotteds{9.75}{11}{4.1}
  \udotteds{9.25}{10.5}{-.1}
  \ldotted{10.25}{11.5}
  \ldotted{-.5-.5+.75}{.75-.25+.75}
  \udotteds{-.5-.5+.75}{.5-.5+.25+.75}{-.1}
  \udotteds{3.625}{4.875}{2.1}
  \draw(4.5,4.1)node[above]{{\small $\overbrace{\phantom{....................}}^{B_{w_x}}$}};
  \draw(9,4.1)node[above]{{\small $\overbrace{\phantom{......................................}}^{\bigcup_{w\in W_2^x}B_w}$}};
  \draw(5.5,1.9)node[below]{{\small $\underbrace{\phantom{..........}}_{D_{u_x}}$}};
  \draw(2.125-.5-.25-.125+.75,0.1)node[below]{{\small $\underbrace{\phantom{............................................}}_{{\bigcup_{v\in V_x^1}A_v}}$}};
  \draw(9,0.1)node[below]{{\small $\underbrace{\phantom{...............................................}}_{{\bigcup_{v\in V_x^2}A_v}}$}};
  \draw(5,6)node[above]{{\small $x$}};
  \draw(8.5,3)node{{\small $\pi_x$}};
  \draw[|-|] (-1-.5,6)--(-1-.5,0);
  \draw[|-|] (-1-.5,4)--(-1-.5,2);
  \draw(-1-.5,5)node[left]{{\small $\al$}};
  \draw(-1-.5,3)node[left]{{\small $\pi$}};
  \draw(-1-.5,1)node[left]{{\small $\be\al$}};
  \muvs504
	\end{tikzpicture}    \caption{A schematic diagram of the product $\al\pi\be\al$, focusing on the transversal block $\{x\}\cup E_x'$, where $\bigcup_{v\in V_x^1}A_v\sub E_x \sub \bigcup_{v\in V_x}A_v$.  See text for further explanation.}
    \label{fig:newnewnewnewprop}
   \end{center}
 \end{figure}
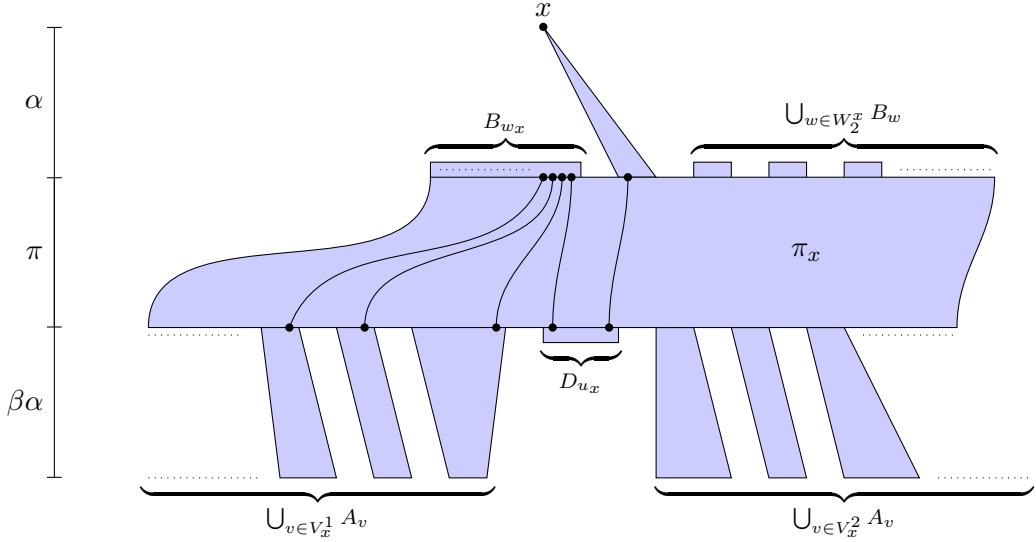

Now fix some $x\in X$.  Consider the bijection
\[
\si_x:\{a_x\}\cup B_{w_x} \to \{b_x,c_x\} \cup \set{e_{xy}}{y\in B_{w_x}\sm\{d_x\}}
\]
defined by $a_x\mt b_x$, $d_x\mt c_x$, and $y\mt e_{xy}$ for each $y\in B_{w_x}\sm\{d_x\}$.  The complements of the domain and codomain of $\si_x$ in $A_x\cup B_{w_x}\cup\bigcup_{w\in W_2^x}B_w$ and $D_{u_x}\cup\bigcup_{v\in V_x}C_v$, respectively, both have cardinality $|X|$, so we may extend $\si_x$ arbitrarily to a bijection
\[
\pi_x:A_x\cup B_{w_x}\cup\bigcup_{w\in W_2^x}B_w \to D_{u_x}\cup\bigcup_{v\in V_x}C_v.
\]
Note that if $\pi\in\S_X$ is any permutation extending $\pi_x$, then the block of $\al\pi\be\al$ containing $x$ is of the form $\{x\}\cup E_x'$ where $E_x\sups\bigcup_{y\in B_{w_x}\sm\{b_x\}}A_{v_{xy}}=\bigcup_{v\in V_x^1}A_v$; see Figure \ref{fig:newnewnewnewprop}.

Now, $\bigcup_{x\in X}\pi_x$ is a bijection from
\[
\bigcup_{x\in X}A_x \cup \bigcup_{x\in X}B_{w_x} \cup \bigcup_{x\in X}\bigcup_{w\in W_2^x}B_w = \bigcup_{x\in X}A_x \cup \bigcup_{z\in Z_{\an}\sm W_2^\infty}B_z
\]
to
\[
\bigcup_{x\in X}D_{u_x}\cup \bigcup_{x\in X}\bigcup_{v\in V_x}C_v = \bigcup_{u\in U_1}D_u\cup\bigcup_{x\in X}C_x.
\]
Since the complements of these sets in $X$ have cardinality $|X|$, we may extend $\bigcup_{x\in X}\pi_x$ arbitrarily to a permutation $\pi\in\S_X$.  Now, $\de=\al\pi\be\al\in\L_X$ satisfies $d^*(\de)\geq d^*(\al)=|X|$, so we may write
$
\de=\LXpartnth.
$
By the above discussion, we see that for each $x\in X$, $|E_x|\geq\left| \bigcup_{v\in V_x^1}A_v\right|\geq|V_x^1|=|B_{w_x}|$.  It follows that, for any cardinal $\mu<|X|$, \begin{align*}
k^*(\de,\mu)&= \# \set{x\in X}{|E_x|\geq\mu}  \geq \# \set{x\in X}{|B_{w_x}|\geq\mu} \\
&= \# \set{w\in W_1}{|B_w|\geq\mu} = |W_1\cap Z_\mu|=|X|,
\end{align*}
as required.  This completes the proof. \epf

\begin{rem}
Some parts of the argument in Case 1 of the above proof are similar to the proof of \cite[Lemma 6.2]{EMP}.
\end{rem}

\smallskip
\begin{thm}\label{thm_singular}
Suppose $|X|$ is singular and let $\al,\be\in\P_X$.  Then $\P_X=\la\S_X,\al,\be\ra$ if and only if (renaming $\al,\be$ if necessary), $\al\in\L_X$, $\be\in\R_X$, $d^*(\al)=d(\be)=|X|$ and either
\bit
	\item[\emph{(i)}] $k^*(\al,2)+d^*(\al,2)=|X|=k(\be,\mu)+d(\be,\mu)$ for all cardinals $\mu<|X|$, or
	\item[\emph{(ii)}] $k^*(\al,\mu)+d^*(\al,\mu)=|X|=k(\be,2)+d(\be,2)$ for all cardinals $\mu<|X|$.
\eit
\end{thm}

\pf The reverse implication was proved in Proposition \ref{newnewnewnewprop}.  So suppose $\P_X=\la\S_X,\al,\be\ra$.  Again, we may assume that $\al\in\L_X$, $\be\in\R_X$, $d^*(\al)=d(\be)=|X|$.  Suppose that (i) and (ii) do not hold.  So one of
\bit
	\item[(I)] $k^*(\al,2)+d^*(\al,2)<|X|$, or
	\item[(II)] $k(\be,\mu)+d(\be,\mu)<|X|$ for some cardinal $\mu<|X|$
\eit
holds, and so too does one of
\bit
	\item[(III)] $k^*(\al,\nu)+d^*(\al,\nu)<|X|$ for some cardinal $\nu<|X|$, or
	\item[(IV)] $k(\be,2)+d(\be,2)<|X|$.
\eit
Again, Lemma \ref{prelim-22} implies that (I) (and, dually, (IV)) cannot hold.  Now suppose (II) and (III) hold.  Let $\lam=\max(\mu,\nu)$.  
Whether $\lam$ is singular or regular, the successor cardinal $\lam^+$ is regular, and we still have $\lam^+<|X|$ as well as
$
 k^*(\al,\lam^+)+d^*(\al,\lam^+)<|X|$ and $k(\be,\lam^+)+d(\be,\lam^+)<|X|.
$
But then Lemma \ref{prelim-nu} implies that $\la\S_X,\al,\be\ra$ does not contain any $\ga\in\L_X$ with $k^*(\ga,\lam^+)=|X|$, a contradiction.  \epf

\section{Relative rank of $\P_X$ modulo $\EX$ and $\EX\cup\S_X$}\label{sect:rankPXEX}

We now turn to the task of calculating the relative rank of $\P_X$ modulo the set of idempotent partitions $\EX=E(\P_X)=\set{\al\in\P_X}{\al=\al^2}$.  We also calculate the relative rank of $\P_X$ modulo the set $\EX\cup\S_X$ of idempotents and units.  We must first recall some ideas from \cite{EF}.  With this in mind, consider a partition
\[
\al = \PXpartn.
\]
We define
\[
s(\al)=\sum_{i\in I}\big(|A_i|-1\big) + \sum_{j\in J}|C_j|
\AND
s^*(\al)=\sum_{i\in I}\big(|B_i|-1\big) + \sum_{k\in K}|D_k|.
\]
These parameters, which were called the \emph{singularity} and \emph{cosingularity} of $\al$ and denoted $\text{sing}(\al)$ and $\text{cosing}(\al)$ in \cite{EF}, allow the alternate characterizations $\L_X=\set{\al\in\P_X}{s(\al)=0}$ and $\R_X=\set{\al\in\P_X}{s^*(\al)=0}$.  Note that $s^*(\al)=s(\al^*)$.  We also write
\[
\sh(\al)=\#\set{i\in I}{A_i\cap B_i=\emptyset}.
\]
This parameter was called the \emph{shift} of $\al$ in \cite{EF}.

For any subset $\Sigma\sub\P_X$, we write $\Sigma^{\fin}$ for the set of all partitions $\al\in\Si$ for which the set $\set{x\in X}{[x]_\al\not=\{x,x'\}}$ is finite.  In \cite{EF}, this set was called the \emph{warp set} of $\al$, and the elements of $\Si^{\fin}$ were called the \emph{finitary} elements of $\Si$.

\bigskip
\begin{thm}[{See \cite[Theorems 30 and 33]{EF}}] \label{EFthm}
Let $X$ be any infinite set.  Then
\bit
\item[\emph{(i)}] $\la \EX\ra = \{1\} \cup \left(\PXfin\sm\SXfin\right) \cup \bigset{\al\in\P_X}{s(\al)=s^*(\al)\geq\max(\an,\sh(\al))}$, and
\item[\emph{(ii)}] $\la \EX\cup\S_X\ra = \bigset{\al\in\P_X}{s(\al)=s^*(\al)}$. \epf
\eit
\end{thm}

\begin{rem}
In the case of finite $X$, we have $\la\EX\cup\S_X\ra=\P_X$ and ${\la\EX\ra=\{1\}\cup(\P_X\sm\S_X)}$.  See \cite[Theorems 32, 36 and 41]{JEgrpm} for presentations of finite $\P_X$ with respect to various generating sets consisting of idempotents and units, including a minimal generating set of size~4.  See \mbox{\cite[Theorem~46]{JEpnsn}} for a presentation of finite $\P_X\sm\S_X$ in terms of a minimal idempotent generating set of size $\frac12|X|\cdot(|X|+1)$.  The minimal generating sets (and minimal idempotent generating sets) for finite $\P_X\sm\S_X$ and various other diagram monoids are classified and enumerated in~\cite{EastGray}.
\end{rem}

In what follows, if $A\sub X$ is any subset, we write
\[
\id_A = \partn{a}{a}{x}{x}_{a\in A,\;\! x\in X\sm A}.
\]

\smallskip
\begin{prop}\label{EXgen}
Let $\al\in\L_X$ and $\be\in\R_X$ be such that $s^*(\al)=s(\be)=|X|$.  Then $\la \EX,\al,\be\ra=\la\EX\cup\S_X,\al,\be\ra=\P_X$.
\end{prop}

\pf Clearly it is sufficient to show that $\la \EX,\al,\be\ra=\P_X$.  Write $\al=\partnL{A_x}{B_i}$ and $\be=\partnR{C_x}{D_j}$.  
First we show that $\la \EX,\al,\be\ra$ contains the symmetric group $\S_X$.  For each $x\in X$, choose some $a_x\in A_x$ and $c_x\in C_x$.  Let $A=\set{a_x}{x\in X}$ and $C=\set{c_x}{x\in X}$, and put
\[
\ga=\al\;\! \id_A = \partnL{a_x}{y}_{x\in X,\;\! y\in X\sm A}
\AND
\de=\id_C\;\! \be = \partnR{c_x}{y}_{x\in X,\;\! y\in X\sm C}.
\]
Note that $|A|=|B|=|X|=|X\sm A|=|X\sm C|$, and that $\ga\ga^*=1=\de^*\de$.  Now let $\pi\in\S_X$ be arbitrary.  Then $\pi=\ga\ga^*\pi\de^*\de$, so it suffices to show that $\ga^*\pi\de^*\in\la\EX\ra$.  Now
\[
\ga^*\pi\de^* = \partn{a_x}{c_{x\pi}}{y}{z}_{x\in X,\;\! y\in X\sm A,\;\! z\in X\sm C}.
\]
So $s(\ga^*\pi\de^*)=|X\sm A|=|X|=|X\sm C|=s^*(\ga^*\pi\de^*)$, and it follows that $\ga^*\pi\de^*\in\la\EX\ra$ by Theorem \ref{EFthm}(i).  This completes the proof that $\S_X\sub\la \EX,\al,\be\ra$.

Now let $(Y_x)_{x\in X\cup\{\infty\}}$ be a moiety of $X\sm A$, where $\infty$ is a symbol that does not belong to $X$.  For each $x\in X$, put $E_x=\{a_x\}\cup Y_x$.  Let
\[
\ve=\partn{E_x}{E_x}{y}{y}_{x\in X,\;\! y\in Y_\infty}.
\]
So $\ve\in\EX$, and
\[
\ga\ve = \partn{x}{E_x}{\emptyset}{y}_{x\in X,\;\! y\in Y_\infty}
\]
belongs to $\L_X$ and satisfies $k^*(\ga\ve,|X|)=|X|=d^*(\ga\ve)$.  Dually, there exists $\eta\in\EX$ such that $\eta\de\in\R_X$ satisfies $k(\eta\de,|X|)=|X|=d(\eta\de)$.  It follows from Proposition \ref{newnewprop} that $\la\EX,\al,\be\ra\sups\la\S_X,\ga\ve,\eta\de\ra=\P_X$. \epf

\begin{lemma}\label{LXEX=RXEX=1}
We have $\L_X\cap\EX=\R_X\cap\EX=\{1\}$.
\end{lemma}

\pf Let $\al\in\L_X\cap\EX$ and write
$\al=\partnL{A_x}{B_i}$.
Then $\al=\al^2$ implies that $A_x=\bigcup_{y\in A_x}A_y$ for all $x\in X$.  This gives $A_x=\{x\}$ for all $x\in X$.  It follows that $I=\emptyset$, and $\al=1$.  A dual argument shows that $\R_X\cap\EX=\{1\}$. \epf

\begin{thm}\label{rankPXEX}
If $X$ is any infinite set, then $\rank(\P_X:\EX)=\rank(\P_X:\EX\cup\S_X)=2$.
\end{thm}

\pf Proposition \ref{EXgen} tells us that $\rank(\P_X:\EX)\leq2$.  Since $\rank(\P_X:\EX\cup\S_X)\leq \rank(\P_X:\S_X)$, it is sufficient to show that $\rank(\P_X:\EX\cup\S_X)\geq2$.  Let $\al\in\P_X$.  The proof will be complete if we can show that $\la\EX\cup\S_X,\al\ra$ is a proper subsemigroup of $\P_X$.  Suppose to the contrary that $\la\EX\cup\S_X,\al\ra=\P_X$.  Let $\be\in\L_X\sm\S_X$, and consider an expression $\be=\ga_1\cdots\ga_r$ where $\ga_1,\ldots,\ga_r\in\EX\cup\S_X\cup\{\al\}$ and $r$ is minimal.  Now $\ga_1\in\L_X$ since $\P_X\sm\L_X$ is a right ideal.  If $\ga_1\in\EX$, we would have $\ga_1=1$ by Lemma \ref{LXEX=RXEX=1}, contradicting either the minimality of $r$ or the fact that $\be\not=1$.  So we must have $\ga_1\in\S_X\cup\{\al\}$.  If $\ga_1\in\S_X$, then $\ga_1^{-1}\be=\ga_2\cdots\ga_r\in\L_X\sm\S_X$, and this expression is also of minimal length.  Continuing in this way, we see that there exists $1\leq s\leq r$ such that $\ga_1,\ldots,\ga_{s-1}\in\S_X$ and $\ga_s=\al$.  So $\al\ga_{s+1}\cdots\ga_r=\ga_{s-1}^{-1}\cdots\ga_1^{-1}\be\in\L_X$, and this implies $\al\in\L_X$.  A dual argument gives $\al\in\R_X$ so that, in fact, $\al\in\S_X$.  But then $\P_X=\la\EX\cup\S_X,\al\ra=\la\EX\cup\S_X\ra$, contradicting Theorem~\ref{EFthm}(ii).~\epf

\begin{rem}
It follows from \cite[Proposition 39]{JEgrpm} and its proof that $\rank(\P_X:\E_X)=2$ if $X$ is finite and $|X|\geq3$.  Since $\P_X=\la\E_X\cup\S_X\ra$ for any finite set $X$ \cite[Theorem 32]{JEgrpm}, it follows that $\rank(\P_X:\E_X\cup\S_X)=0$ for finite $X$.  
\end{rem}

\section{Generating pairs for $\P_X$ modulo $\EX$ and $\EX\cup\S_X$}\label{sect:genpairsPXEX}

In order to establish the converse of Proposition \ref{EXgen}, we will first need to prove a series of lemmas.

\bigskip
\begin{lemma}\label{EXalbelem}
Suppose $\al,\be\in\P_X$ are such that $\P_X=\la\EX\cup\S_X,\al,\be\ra$.  Then (renaming $\al,\be$ if necessary) $\al\in\L_X$ and $\be\in\R_X$.
\end{lemma}

\pf A similar argument to that in the proof of Theorem \ref{rankPXEX} shows that one of $\al,\be$, say~$\al$, belongs to $\L_X$, and a dual argument shows that one of $\al,\be$ belongs to $\R_X$.  If, in fact, $\al\in\R_X$ as well, then $\al\in\S_X$, and we would have $\P_X=\la\EX\cup\S_X,\al,\be\ra=\la\EX\cup\S_X,\be\ra$, contradicting Theorem \ref{rankPXEX}.  So it follows that $\be\in\R_X$. \epf

A \emph{quotient} of $X$ is a collection $\bY=\set{A_i}{i\in I}$ of pairwise disjoint nonempty subsets of $X$ such that $X=\bigcup_{i\in I}A_i$.  We write $\bY\pre X$ to indicate that $\bY$ is a quotient of $X$.  If $\bY\pre X$ is as above, we write
\[
\id_\bY = \partn{A_i}{A_i}{\emptyset}{\emptyset}_{i\in I}.
\]

\smallskip
\begin{lemma}[{See \cite[Lemma 5 and Proposition 6]{EF}}]\label{IXJX}
If $X$ is any infinite set, then $\la\EX\ra=\la\Si\ra,$ where
$\Si = \set{\id_A}{A\sub X} \cup \set{\id_\bY}{\bY\pre X}$. \epfres
\end{lemma}

The previous result will substantially simplify the proof of the following technical lemma.

\bigskip
\begin{lemma}\label{prelim-s*}
Suppose $\al\in\L_X$ and $\be\in\R_X$ and that $s^*(\al)<|X|$.  If $\ga_1,\ldots,\ga_r\in\EX\cup\S_X\cup\{\al,\be\}$ are such that $\ga_1\cdots\ga_r\in\L_X$, then $s^*(\ga_1\cdots\ga_r)<|X|$.
\end{lemma}

\pf Write $\al=\partnL{A_x}{B_i}$ and $\be=\partnR{C_x}{D_j}$.  
The $r=1$ case is trivial, so suppose $r\geq2$ and put $\ga=\ga_1\cdots\ga_{r-1}$.  Since $\ga\in\L_X$, an inductive hypothesis gives $s^*(\ga)<|X|$.  Write
$\ga=\partnL{E_x}{F_k}$.
Since $s^*(\al)<|X|$ and $s^*(\ga)<|X|$, it follows that $\sum_{x\in X}\big(|A_x|-1\big)$, $\sum_{i\in I}|B_i|$, $\sum_{x\in X}\big(|E_x|-1\big)$ and $\sum_{k\in K}|F_k|$ are all less than $|X|$.  We now consider four separate cases according to whether $\ga_r\in\EX$, $\ga_r\in\S_X$, $\ga=\al$ or $\ga=\be$.  In each case, we must show that $s^*(\ga\ga_r)<|X|$.

{\bf Case 1.}  First suppose $\ga_r\in\EX$.  In fact, by Lemma \ref{IXJX}, we may assume that $\ga_r=\id_A$ for some $A\sub X$ or $\ga_r=\id_\bY$ for some $\bY\pre X$.

{\bf Subcase 1.1.}  Suppose $\ga_r=\id_A$ for some $A\sub X$.  For $x\in X$ and $k\in K$, let $G_x=E_x\cap A$ and $H_k=F_k\cap A$.  Then
\[
\ga\ga_r = \partnL{G_x}{H_k,y}_{x\in X,\;\! k\in K,\;\! y\in X\sm A}.
\]
Note that some of the $H_k$ may be empty, but all of the $G_x$ are nonempty.  Now
\begin{align*}
s^*(\ga\ga_r) &= \sum_{x\in X}\big(|G_x|-1\big) + \sum_{k\in K}|H_k| +|X\sm A| \\
&= \sum_{x\in X} \big( |G_x|-1 + |(X\sm A)\cap E_x|\big) + \sum_{k\in K} \big(|H_k|+|(X\sm A)\cap F_k|\big) \\
&= \sum_{x\in X}\big(|E_x|-1\big) + \sum_{k\in K}|F_k| = s^*(\ga)<|X|.
\end{align*}

{\bf Subcase 1.2.}  Suppose $\ga_r=\id_\bY$ for some $\bY\pre X$.  Let $\ve$ be the equivalence relation on $X$ corresponding to $\bY$.  That is, two elements of $X$ are $\ve$-related if and only if they belong to the same block of $\bY$, and we have $\bY=X/\ve$.  Put $\eta=\coker(\ga)$, and let $\bZ=X/\eta$.  Clearly we have $\ga=\ga\;\!\id_\bZ$.  So $\ga\ga_r=\ga\;\!\id_\bZ\;\!\id_\bY=\ga\;\!\id_\bW$, where $\bW=X/(\ve\vee\eta)$; here $\ve\vee\eta$ denotes the least equivalence on $X$ containing $\ve\cup\eta$.  Since $\eta\sub\ve\vee\eta$, every block of $\bZ$ is contained in a block of $\bW$.  Since $\ga\ga_r\in\L_X$, it is not possible for a block of $\bW$ to contain $E_{x_1}$ and $E_{x_2}$ if $x_1\not=x_2$.  So we may write $\bW=\set{U_x}{x\in X}\cup\set{V_l}{l\in L}$ where $E_x\sub U_x$ for all $x\in X$, and we have
$\ga\ga_r = \partnL{U_x}{V_l}$.
Now, for each $x\in X$, there exists a subset $K_x\sub K$ such that $U_x=E_x\cup\bigcup_{k\in K_x}F_k$.  And for each $l\in L$, there exists a subset $K_l\sub K$ such that $V_l=\bigcup_{k\in K_l}F_k$.  Note that $K=\bigcup_{x\in X}K_x\cup\bigcup_{l\in L}K_l$.  Then
\begin{align*}
s^*(\ga\ga_r) = \sum_{x\in X}\big(|U_x|-1\big) + \sum_{l\in L}|V_l| 
&= \sum_{x\in X} \left( |E_x|-1 + \sum_{k\in K_x}|F_k| \right) + \sum_{l\in L}\sum_{k\in K_l}|F_k| \\
&= \sum_{x\in X}\big(|E_x|-1\big) +\sum_{k\in K}|F_k|=s^*(\ga)<|X|.
\end{align*}

{\bf Case 2.}  If $\ga_r\in\S_X$, then clearly $s^*(\ga\ga_r)=s^*(\ga)<|X|$.

{\bf Case 3.}  Next suppose $\ga_r=\al$.  Now $\ga\al = \partnL{G_x}{B_i,H_k}$, where
where $G_x=\bigcup_{y\in E_x}A_y$ and $H_k=\bigcup_{y\in F_k}A_y$ for each $x,k$.  Let $Y=\set{x\in X}{|A_x|\geq2}$.  Since ${\sum_{x\in X}\big(|A_x|-1\big)<|X|}$, it follows that $|Y|<|X|$.  So ${\sum_{x\in X}|A_x|=\sum_{x\in X}\big(|A_x|-1\big)+|Y|<|X|}$.  But then
\begin{align*}
s^*(\ga\al) &= \sum_{x\in X}\big(|G_x|-1\big) +\sum_{k\in K}|H_k| +\sum_{i\in I}|B_i| \\
&= \sum_{x\in X}\left( \sum_{y\in E_x}|A_y| - 1 \right) + \sum_{k\in K}\sum_{y\in F_k}|A_y| + \sum_{i\in I}|B_i| \\
&\leq \sum_{x\in X}\sum_{y\in E_x}|A_y| + \sum_{k\in K}\sum_{y\in F_k}|A_y| + \sum_{i\in I}|B_i| \\
&= \sum_{x\in X}|A_x| + \sum_{i\in I}|B_i| <|X|.
\end{align*}

{\bf Case 4.}  Finally, suppose $\ga_r=\be$.  Write $\ga\be=\partnL{P_x}{Q_l}$.
Let $x\in X$.  The middle row of the connected component containing $x$ in the product graph $\Ga(\ga,\be)$ is (omitting double dashes)
\[
E_x\cup\bigcup_{k\in K_x}F_k=\bigcup_{y\in P_x}C_y\cup\bigcup_{j\in J_x}D_j
\]
for some subsets $K_x\sub K$ and $J_x\sub J$.  So $|P_x|\leq|E_x|+\sum_{k\in K_x}|F_k|$.  

Now let $l\in L$.  The middle row of the connected component containing $Q_l'$ in the product graph $\Ga(\ga,\be)$ is (omitting double dashes)
\[
\bigcup_{k\in K_l}F_k = \bigcup_{y\in Q_l}C_y \cup \bigcup_{j\in J_l}D_j
\]
for some subsets $K_l\sub K$ and $J_l\sub J$.  Thus, $|Q_l|\leq\sum_{k\in K_l}|F_k|$.  It follows that
\begin{align*}
s^*(\ga\be) = \sum_{x\in X}\big(|P_x|-1\big) + \sum_{l\in L}|Q_l| 
&\leq \sum_{x\in X}\left( |E_x|-1+\sum_{k\in K_x}|F_k| \right) + \sum_{l\in L}\sum_{k\in K_l}|F_k| \\
&\leq \sum_{x\in X}\big(|E_x|-1\big)+\sum_{k\in K}|F_k|=s^*(\ga)<|X|.
\end{align*}
This completes the proof. \epf

\begin{rem}
Some elements of the argument from Subcase 1.2 are similar to the proof of \cite[Lemma 14]{EF}.
\end{rem}

\smallskip
\begin{thm}\label{EXgenthm}
Let $X$ be any infinite set and let $\al,\be\in\P_X$.  Then $\P_X=\la\EX,\al,\be\ra$ if and only if $\P_X=\la\EX\cup\S_X,\al,\be\ra$ if and only if (renaming $\al,\be$ if necessary) $\al\in\L_X$ and $\be\in\R_X$ satisfy $s^*(\al)=s(\be)=|X|$.
\end{thm}

\pf In Proposition \ref{EXgen}, we saw that if $\al\in\L_X$ and $\be\in\R_X$ satisfy ${s^*(\al)=s(\be)=|X|}$, then $\P_X=\la\EX,\al,\be\ra$.  It is obvious that $\P_X=\la\EX,\al,\be\ra$ implies $\P_X=\la\EX\cup\S_X,\al,\be\ra$.  Suppose now that $\P_X=\la\EX\cup\S_X,\al,\be\ra$.  By Lemma \ref{EXalbelem}, we may assume that $\al\in\L_X$ and $\be\in\R_X$.  If $s^*(\al)<|X|$, then Lemma \ref{prelim-s*} would imply that any element $\ga\in\L_X$ with $s^*(\ga)=|X|$ could not belong to $\la\EX\cup\S_X,\al,\be\ra$, a contradiction.  Thus, $s^*(\al)=|X|$.  A dual argument shows that $s(\be)=|X|$.  This completes the proof. \epf

\section{Sierpi\'nski rank and the semigroup Bergman property}\label{sect:sierpinski}

Let $S$ be a semigroup.  Recall that the \emph{Sierpi\'nski rank} of $S$, denoted $\SR(S)$, is the least integer~$n$ such that every countable subset of $S$ is contained in an $n$-generator subsemigroup of $S$, if such an integer exists.  Otherwise, we say that $S$ has infinite Sierpi\'nski rank. 

\bigskip
\begin{thm}\label{PXsierpinski}
Let $X$ be any infinite set.  Then $\SR(\P_X)\leq4$.
\end{thm}

\pf A general result \cite[Lemma 2.3]{MPsurj} states that if $T$ is a subsemigroup of a semigroup~$S$, then $\SR(S)\leq\rank(S:T)+\SR(T)$ if $\rank(S:T)$ and $\SR(T)$ are finite.  By Theorem \ref{cor4}, and the fact that $\SR(\S_X)=2$ \cite[Theorem 3.5]{Gal}, it immediately follows that $\SR(\P_X)\leq4$.  But for the sake of completeness, and since it will be useful in a subsequent proof, we offer a direct proof that is reminiscent of Banach's proof \cite{Ban} that the full transformation semigroup~$\T_X$ has Sierpi\'nski rank $2$ \cite{Sie}, and is also similar to the proof of \cite[Proposition 4.2]{HHMR}.

\begin{figure}[ht]
   \begin{center}
\begin{tikzpicture}[xscale=1.2,yscale=0.4]
  \bpermblock0224
  \bpermblock2446
  \bpermblock4668
  \ldotted{8.25}{10}
  \draw(0,0)node{\LARGE {$.$}};
  \draw(.125,0)node{\LARGE {$.$}};
  \draw(.25,0)node{\LARGE {$.$}};
  \ldotted{.5}{1.75}
  \udotteds{6.25}{10}4
  \draw(1,6)node[above]{{\small $X_0$}};
  \draw(3,6)node[above]{{\small $X_1$}};
  \draw(5,6)node[above]{{\small $X_2$}};
  \draw(2,3)node{{\small $\phi_1$}};
  \draw(4,3)node{{\small $\phi_2$}};
  \draw(6,3)node{{\small $\phi_3$}};
	\draw(-0.2,3)node[left]{$\be=$}; 
	 \draw(7,-.5)node[below]{{\tiny $\phantom{|X|}$}}; 
	\end{tikzpicture}
	~ ~ ~
	\begin{tikzpicture}[xscale=1.2,yscale=0.4]
	\draw(-0.2,3)node[left]{$\ga=$}; 
  \bpermblock240{.4}
  \bpermblock46{.4}{.8}
  \bpermblock68{.8}{1.2}
  \bpermblock0{.4}24
  \bpermblock{.4}{.8}46
  \bpermblock{.8}{1.2}68
  \bpermblockd240{.4}
  \bpermblockd68{.8}{1.2}
  \ldotted{8.25}{10}
  \ldotted{1.45}{1.75}
  \udotteds{8.25}{10}4
  \udotteds{1.45}{1.75}4
  \draw(3,6)node[above]{{\small $X_1$}};
  \draw(5,6)node[above]{{\small $X_2$}};
  \draw(7,6)node[above]{{\small $X_3$}};
  \draw(.2,6)node[above]{{\small $Y_1$}};
  \draw(.6,6)node[above]{{\small $Y_2$}};
  \draw(1,6)node[above]{{\small $Y_3$}};
  \draw(3-.2,1)node{{\small $\de_1$}};
  \draw(5-.2,1)node{{\small $\de_2$}};
  \draw(7-.2,1)node{{\small $\de_3$}};
  \draw(3-.2,5)node{{\small $\psi_1$}};
  \draw(5-.2,5)node{{\small $\psi_2$}};
  \draw(7-.2,5)node{{\small $\psi_3$}};
	\end{tikzpicture}
    \caption{The partitions $\be$ (top) and $\ga$ (bottom) from the proof of Theorem~\ref{PXsierpinski}.}
    \label{fig:sierpinski}
   \end{center}
 \end{figure}
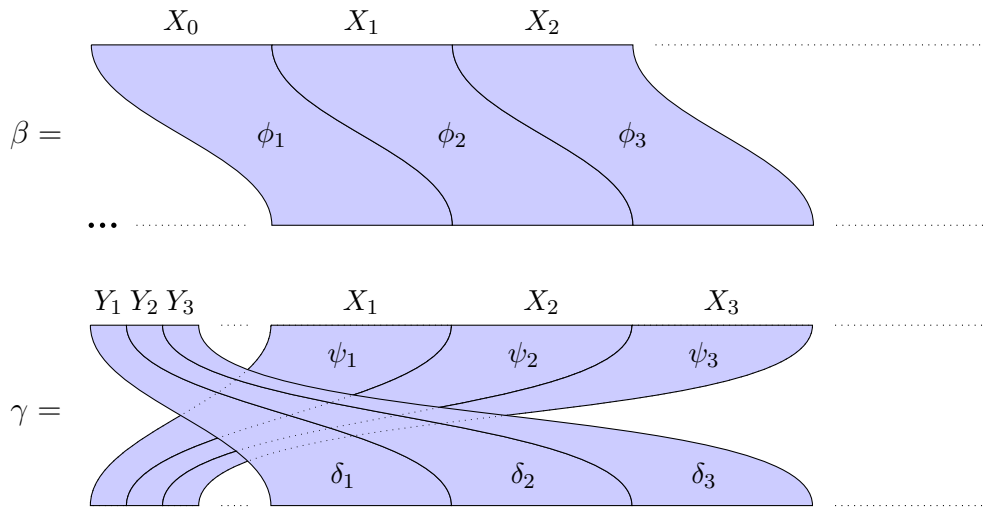

With this in mind, suppose we have a countable subset $\Si=\set{\al_n}{n\in\N}$ of $\P_X$.  For $n\in\N$, write
\[
\al_n = \partn{A_i^n}{B_i^n}{C_j^n}{D_k^n}_{i\in I^n,\;\! j\in J^n,\;\! k\in K^n}.
\]
We will construct two partitions $\be,\ga\in\P_X$ such that $\Si\sub\la\be,\be^*,\ga,\ga^*\ra$.  
Let $(X_n)_{n\in\N\cup\{0\}}$ be a moiety of $X$, and let $(Y_n)_{n\in\N}$ be a moiety of $X_0$.  For each $n\in\N$, fix bijections $\phi_n:X_{n-1}\to X_n$ and $\psi_n:X_n\to Y_n$.  Let $\phi=\bigcup_{n\in\N}\phi_n$ and $\psi=\bigcup_{n\in\N}\psi_n$, noting that these are bijections $\phi:X\to X\sm X_0$ and $\psi:X\sm X_0\to X_0$.  For each $n\in\N$, define $\si_n=\phi\psi\phi^n$ and $\tau_n=\phi\psi\phi^n\psi$, noting that these are bijections $\si_n:X\to X_n$ and $\tau_n:X\to Y_n$.  For each $n\in\N$, also define
\[
\de_n=\partn{A_i^n\tau_n}{B_i^n\si_n}{C_j^n\tau_n}{D_k^n\si_n}_{i\in I^n,\;\! j\in J^n,\;\! k\in K^n}.
\]
(Note that $\de_n$ is not a full partition.  Rather, $\bigcup_{i\in I^n}A_i^n\tau_n \cup \bigcup_{j\in J^n}C_j^n\tau_n = Y_n$ and $\bigcup_{i\in I^n}B_i^n\tau_n \cup \bigcup_{k\in K^n}D_k^n\tau_n = X_n$.)  Now put
\[
\be = \partn{x}{x\phi}{\emptyset}{y}_{x\in X,\;\! y\in X_0} \AND \ga=\bigcup_{n\in\N}\psi_n\cup\bigcup_{n\in\N}\de_n.
\]
See Figure \ref{fig:sierpinski}.  One may easily check that $\be\ga\be^n\ga^2(\be^*)^n\ga^*\be^*=\al_n$ for each $n\in\N$. \epf

\begin{cor}\label{PXsierpinski_cor}
If $\Si$ is any subset of $\P_X$, then $\rank(\P_X:\Si)$ is either uncountable or at most~$4$.
\end{cor}

\pf If $\rank(\P_X:\Si)\leq\an$, then $\P_X=\la\S_X\cup\Gamma\ra$ for some countable subset $\Ga\sub\P_X$.  But, by Theorem \ref{PXsierpinski}, $\Gamma\sub\la\Lambda\ra$ for some $\Lambda\sub\P_X$ with $|\Lambda|\leq4$.  But then $\P_X=\la\Si\cup\Lambda\ra$. \epf

\begin{rem}
A recent result of Hyde and P\'eresse \cite[Theorem 1.4]{HP} shows that the Sierpi\'nski rank of an infinite symmetric inverse monoid is equal to $2$, an improvement of {\cite[Proposition 4.2]{HHMR}} which gave an upper bound of $4$.  It is anticipated that the methods of \cite{HP} may be extended to show that $\SR(\P_X)=2$, but this is beyond the scope of the current work.  Naturally, this would show that Corollary \ref{PXsierpinski_cor} could be suitably improved too.
\end{rem}


Recall that a semigroup $S$ has the \emph{semigroup Bergman property} \cite{MMR} if the length function of $S$ is bounded with respect to any generating set for $S$.  The property has this name since Bergman showed in \cite{Berg} that an infinite symmetric group $\S_X$ has the property.  (Actually, Bergman showed that $\S_X$ has this property with respect to \emph{group} generating sets of $\S_X$, and the semigroup analogue was shown in \cite[Corollary 2.5]{MMR}.)

Recall from \cite{MMR} that a semigroup $S$ is said to be \emph{strongly distorted} if there exists a sequence of natural numbers $(a_n)_{n\in\N}$, and a natural number $N_S$ such that, for all sequences $(s_n)_{n\in\N}$ of elements of $S$, there exist $t_1,\ldots,t_{N_S}\in S$ such that each $s_n$ can be expressed as a product of length at most $a_n$ in the elements $t_1,\ldots,t_{N_S}$.

\bigskip
\begin{prop}[{See \cite[Lemma 2.4 and Proposition 2.2(i)]{MMR}}]\label{MMRprop}
If $S$ is non-finitely generated and strongly distorted, then $S$ has the semigroup Bergman property. \epfres
\end{prop}

Since $|\P_X|>\an$ for any infinite set $X$, $\P_X$ is clearly not finitely generated.  And the proof of Theorem \ref{PXsierpinski} shows that $\P_X$ is strongly distorted (we take $N_{\P_X}=4$ and $a_n=2n+6$ for all $n\in\N$).  So we immediately obtain the following.

\bigskip
\begin{thm}\label{PXbergman}
If $X$ is any infinite set, then $\P_X$ has the semigroup Bergman property. \epfres
\end{thm}

\footnotesize
\def\bibspacing{-1.1pt}
\bibliography{ipms}
\bibliographystyle{plain}
\end{document}